\documentclass{amsart}
\usepackage[latin1]{inputenc}
\usepackage{color}
\usepackage{amssymb}
\usepackage{enumerate}
\usepackage[all]{xy}

\usepackage{mathrsfs}
\usepackage{bm}
\usepackage{bbm}
\usepackage{enumitem}

\usepackage{hyperref,cleveref,graphics,mathrsfs}
\newtheorem{thm}{Theorem}[section]
\newtheorem{cor}[thm]{Corollary}
\newtheorem{lem}[thm]{Lemma}
\newtheorem{prop}[thm]{Proposition}
\theoremstyle{definition}

\newtheorem{rem}[thm]{Remark}

\numberwithin{equation}{section}
\newcommand{\Real}{\mathbb R}
\newcommand{\FF}{\mathcal F}
\newcommand{\CC}{\mathbf C}
\newcommand{\CO}{\overline{\mathbf C}}

\newcommand{\Nat}{\mathbb{N}}

\def\epsilon{\varepsilon}

\newcommand{\timur}[1]{{\textcolor{blue}{{\bf T.O:} #1}}}
\newcommand{\marg}[1]{\marginpar{\tiny #1}}     % A margin note
\usepackage{bbm}
\newcommand{\vertiii}[1]{{\left\vert\kern-0.25ex\left\vert\kern-0.25ex\left\vert #1 
    \right\vert\kern-0.25ex\right\vert\kern-0.25ex\right\vert}}

\newcommand{\triple}[1]{{\left\vert\kern-0.25ex\left\vert\kern-0.25ex\left\vert #1 
    \right\vert\kern-0.25ex\right\vert\kern-0.25ex\right\vert}}

\newcommand{\sphere}{\mathbf{S}}
\newcommand{\ball}{\mathbf{B}}

\newcommand{\vr}{\varepsilon}

\newcommand{\dist}{{\mathrm{dist}}}
\newcommand{\sgn}{{\mathrm{sign}}}
\newcommand{\supp}{{\mathrm{supp}}}
\newcommand{\diam}{{\mathrm{diam}} \, }

\newcommand{\spn}{\mathrm{span}}
\newcommand{\fbl}{\mathrm{FBL}}
\newcommand{\fvl}{\mathrm{FVL}}
\newcommand{\fbp}{{\mathrm{FBL}}^{(p)}}
\newcommand{\pfbp}{{\mathbf{X}}_p}
\newcommand{\css}{{\mathcal{C}}}
\newcommand{\boldom}{\boldsymbol{\omega}}

\newcommand{\ran}{\textrm{ran }}

\title[Unit balls of free Banach lattices]{Geometry of unit balls of free Banach lattices, and its applications}



\author{T.~Oikhberg}\address{Dept. of Mathematics, University of Illinois, Urbana IL 61801, USA} \email{oikhberg@illinois.edu}

\date{\today}

\thanks{The author is grateful to the organizers of Workshop on Banach spaces and Banach lattices II at ICMAT (May 2022), where part of this work was carried out. Many interesting conversations with A.~Chavez-Dominguez, P.~Tradacete, and V.~Troitsky were also vital. Last but not least, the author wishes to sincerely thank the anonymous referee for their helpful remarks and suggestions.}

\begin{document}

\begin{abstract}
We begin by describing the unit ball of the free $p$-convex Banach lattice over a Banach space $E$ (denoted by $\fbp[E]$) as a closed solid convex hull of an appropriate set. Based on it, we show that, if a Banach space $E$ has the $\lambda$-Approximation Property, then $\fbp[E]$ has the $\lambda$-Positive Approximation Property. Further, we show that operators $u \in B(E,F)$ (where $E$ and $F$ are Banach spaces) which extend to lattice homomorphisms from $\fbl^{(q)}[E]$ to $\fbp[F]$ are precisely those whose adjoints are $(q,p)$-mixing. Related results are also obtained for free lattices with an upper $p$-estimate.
\end{abstract}

\date{\today}
\maketitle

\tableofcontents

\section{Introduction and preliminaries}\label{intro}

In recent years, the theory of free Banach lattices has been attracting a lot of attention. To recall the definitions, fix a Banach space $E$, and $p \in [1,\infty]$. We say that a Banach lattice $Z$ is the \emph{free $p$-convex Banach lattice over $E$} if
\begin{enumerate}
 \item $Z$ is $p$-convex with constant $1$ -- that is, for any $z_1, \ldots, z_n \in Z$, we have $\| (\sum_i |z_i|^p \big)^{1/p} \| \leq ( \sum_i \|z_i\|^p )^{1/p}$.
 \item There exists a linear isometry $\phi_E : E \to Z$, so that $\phi_E(E)$ generates $Z$ as a Banach lattice.
 \item For any Banach lattice $X$, which is $p$-convex with constant $1$, and any contraction $T : E \to X$, there exists a (necessarily unique) contractive lattice homomorphism $\widehat{T} : Z \to X$ so that $\widehat{T} \phi_{E,p} = T$ -- in other words, the diagram below commutes.
$$
\xymatrix{Z\ar[rd]^{\widehat{T}}&\\
     E\ar[r]^{T}\ar^{\phi_{E,p}}[u]& X}
$$
\end{enumerate}
We denote this free lattice by $\fbp[E]$, and use $\fbl[E]$ for $\fbl^{(1)}[E]$. We simply write $\phi$ for $\phi_{E,p}$ when no confusion can arise.

If $p \in (1,\infty)$, we similarly say that $Z$ is the \emph{free  upper $p$-estimate Banach lattice over $E$} if similar conditions are satisfied, with $\phi_{E, \uparrow p}$ taking the place of $\phi_{E,p}$, and $p$-convexity replaced by having upper $p$-estimate with constant $1$ -- that is, for any $z_1, \ldots, z_n \in Z$, we have $\| \vee_i |z_i| \| \leq ( \sum_i \|z_i\|^p )^{1/p}$. We denote such a lattice by $\fbl^{\uparrow p}[E]$; again, $\phi_{E, \uparrow p}$ is often replaced by $\phi$.

In addition to being interesting in their own right, free Banach lattices have been essential in resolving a question of J.~Diestel: \cite{ART} provides an example of a free Banach lattice which is weakly compactly generated as a Banach lattice, but not in the Banach space sense.
We refer the reader to \cite{JLTTT} and \cite{OiTTT1} for more information about free lattices.

In this paper, we begin by examining, in \Cref{s:fbp hull}, the geometry of the unit ball of a free Banach lattice. The relevant notions (of a solid convex hull etc.) are recalled in \Cref{ss:hulls}. In \Cref{ss:fbp hull} (\Cref{hull}) we represent the unit ball of $\fbp[E]$ as a closed solid convex hull of a certain ``manageable'' set. In \Cref{ss:fb upper p hull} (\Cref{hull upper est}) we examine the same problem for $\fbl^{\uparrow p}[E]$.

The geometric information described above is then applied to several problems, starting with approximation properties of free Banach lattices (\Cref{bap}). Recall that a Banach space $E$ has the \emph{$\lambda$-Approximation Property} (\emph{$\lambda$-AP} for short) if there exists a net of finite rank operators $(u_\alpha) \subset B(E)$ so that $\sup \|u_\alpha\| \leq \lambda$, and $u_\alpha \to I_E$ point-norm. Equivalently, for every finite dimensional $G \subset E$, and any $\vr > 0$, there exists a finite rank $u \in B(E)$ so that $u|_G = I_G$, and $\|u\| < \lambda + \vr$. % The term \emph{Metric Approximation Theory} (\emph{MAP}) is used to denote 1-AP.
If $E$ is a Banach lattice, and the operators $u_\alpha$ (or $u$) as above are positive, then we can talk about the \emph{$\lambda$-Positive Approximation Property} (\emph{$\lambda$-PAP}). %, or the \emph{Positive Metric Approximation Theory} (\emph{PMAP}).

The main result of the section (\Cref{t:bap to pbap}) states that, if $E$ has the $\lambda$-AP, then $\fbp[E]$ has the $\lambda$-PAP. From this, we conclude that, if $\fbp[E]$ has the $\lambda$-AP, then it also has the $\lambda$-PAP \Cref{c:bap to pbap}; similar results also hold for $\fbl^{\uparrow p}[E]$.
In general, it is an open problem whether a Banach lattice possessing the $\lambda$-AP must also have the $\nu$-PAP, for some $\nu$.

\Cref{s:duality} examines the dual lattices $\fbp[E]^*$ and $\fbl^{\uparrow p}[E]^*$. Specifically, \Cref{fbp dual} gives an explicit formula for linear combinations of atoms in $\fbp[E]^*$. The relevant results about $\fbl^{\uparrow p}[E]^*$ also hold (\Cref{fbl p-conv dual}).

In \Cref{s:embeddings} we ``localize'' norm computations in $\fbp[E]$ (\Cref{p:evaluate norm}).

Finally, in \Cref{s:maps} we examine ``how different'' the lattices $\fbp[E]$ can be, for different values of $p$.
\Cref{mixing} shows that an operator $u : E \to F$ $u$ extends to a bounded lattice homomorphism $\overline{u} : \fbl^{(q)}[E] \to \fbp[F]$ if and only if $u^*$ is $(q,p)$-mixing. Consequently (\Cref{t:fin dim formal id}), $\|id : \fbl^{(q)}[E] \to \fbp[E]\| \leq n^\alpha$, where $$\alpha = \alpha(p,q) = \max\Big\{0, \max\Big\{1,\frac{p}2\Big\}\Big(\frac1p - \frac1q\Big)\Big\};$$ additionally, by \Cref{fbp to upper p}, $\|id : \fbp[E] \to \fbl^{\uparrow p}[E]\| \prec (\ln n)^{1/p}$. This information can be used to estimate the ``lattice Banach-Mazur distance'' between free lattices (Remarks \ref{r:Banach Mazur} and \ref{r:optimality of log}).

Throughout this paper, we freely use the standard functional analysis facts and notation. For information, we refer the reader to \cite{AB} or \cite{M-N} (Banach lattices), as well as \cite{alb-kal}, \cite{LT1}, or \cite{LT2} (Banach spaces; the last cited monograph deals with Banach lattices as well).

There are also some notations specific to this paper. The closed unit ball (the unit sphere) of $E$ is denoted by $\ball(E)$ (resp.~$\sphere(E)$). If $e \in E$, we shall sometimes write $\delta_e$ for its image $\phi(e)$ in a free Banach lattice over $E$ (this notation is common in the theory of free lattices). For indices, we use $p' = p/(p-1)$ (so $1/p + 1/p' = 1$). If $A$ and $B$ are expressions in certain variables, we write $A \prec B$ if there is a constant $C$ (perhaps depending on the parameter $p$, but on nothing else) so that $A \leq CB$ no matter what the values of the variables are. $A \sim B$ means that both $A \prec B$ and $A \succ B$ hold.

Concluding this introduction, we mention a functional representation of $\fbp[E]$. Denote by $H[E]$ the space of positively homogeneous functions on $E^*$, which are weak$^*$-continuous on bounded sets. Define $H[E] \ni \delta_e : e^* \mapsto \langle e^*, e \rangle$ and $\phi : E \to H[E] : e \mapsto \delta_e$; $\fbp[E]$ is the closure of $\phi(E) = \{ \delta_e : e \in E\}$ in $H[E]$, equipped with the norm
$$
\|f\| := \sup \Big\{ \big( \sum_{i=1}^N |f(e_i^*)|^p \big)^{1/p} : \|(e_i^*)_{i=1}^N\|_{p,{\textrm{weak}}} \leq 1 \Big\} ,
$$
where
$$
\|(e_i^*)_{i=1}^N\|_{p,{\textrm{weak}}} = \sup \Big\{ \big( \sum_{i=1}^N | \langle e_i^* , e \rangle|^p \big)^{1/p} : e \in E, \|e\| \leq 1 \Big\}
$$
(see \cite{OiTTT1} for more on this).

\section{Representing unit balls as solid convex hulls}\label{s:fbp hull}

\subsection{Survey of solid convex hulls}\label{ss:hulls}

% In this section, we describe the unit balls of free lattices as solid closed convex hulls of ``manageable'' sets. 
% For a set $A$ in a normed space $Z$, we denote by ${\mathrm{CH}}(A)$ its convex hull. If, in addition, $Z$ is a vector lattice, we denote by ${\mathrm{S}}(A)$ its \emph{solid hull} -- that is, the set of all $z \in Z$ for which there exists $a \in A$ with $|z| \leq |a|$.
Recall that a subset $A$ of a vector lattice $Z$ is called \emph{solid} if $z \in A$ whenever $|z| \leq |a|$ for some $a \in A$. Following \cite{OiTu}, we define \emph{the closed solid convex hull} (denoted by ${\mathrm{CSCH}}(A)$) of $A$ as the smallest closed solid convex subset of $Z$ containing $A$. We refer the reader to \cite{OiTu} for more information about such hulls. For future use, we need some basic facts. The first lemma may be known to experts, but we have not seen it stated explicitly.

\begin{lem}\label{solid closed}
 If $B$ is a solid convex subset of $Z$, then so is its norm closure $\overline{B}$.
\end{lem}

\begin{proof}
 The convexity of the closure is well-known. To show the solidity, note first that, if $b_i \to b$, then $|b_i| \to |b|$ as well (the ``triangle inequality''), and therefore, if $b \in \overline{B}$, then $|b| \in \overline{B}$ as well.
 
 It remains to establish that, if $b \in \overline{B}$ is non-negative, then every $x \in [-b,b]$ belongs to $\overline{B}$ as well.
 To this end, write $x = x_+ - x_-$, with $x_+, x_- \in [0,b]$. For $\varepsilon > 0$, and find $c \in B$ with $\|b-c\| < \varepsilon/2$. As $B$ is solid, we can assume (by replacing $c$ by $c_+$) that $c \geq 0$.
 Now let $y_+ = x_+ \wedge c$ and $y_- = x_- \wedge c$. As $x_+ = x_+ \wedge b$, we have $\|x_+ - y_+\| \leq \|b-c\| < \varepsilon/2$, and likewise, $\|x_- - y_-\| < \varepsilon/2$. Clearly, $y_+ - y_- \in [-c,c] \subset B$, and $\|x - y\| < \varepsilon$.
\end{proof}

\begin{lem}\label{separate}
Suppose $A$ is a bounded subset of $Z_+$, where $Z$ a Banach lattice. Then, for $z^* \in Z^*_+$, we have
$$
\sup z^* \big|_{{\mathrm{CSCH}}(A)} = \sup z^*|_A .
$$
\end{lem}

\begin{proof}
The set $\displaystyle \big\{ x \in Z : | \langle z^* , x \rangle | \leq \sup z^*|_A \big\}$
is clearly solid, convex, closed, and contains $A$, hence it also contains ${\mathrm{CSCH}}(A)$. Thus,
$\displaystyle \sup z^* \big|_{{\mathrm{CSCH}}(A)} \leq \sup z^*|_A$; the converse inequality is straightforward.
%
%  Combining \Cref{solid closed} with \cite[Proposition 19.4]{OiTu}, we conclude that ${\mathrm{CSCH}}(A)$ is the closure of the solid hull of the convex hull of $A$, the latter being denoted by ${\mathrm{CH}}(A)$. Therefore,
%  $$ \sup z^* \big|_{{\mathrm{CSCH}}(A)} = \sup z^* \big|_{{\mathrm{CH}}(A)} , $$
% and the latter is clearly equal to $\sup z^*|_A$.
\end{proof}

We also need a ``positive separation'' result.

\begin{lem}\label{positive separation}
 Suppose $A$ is a closed solid convex subset of a Banach lattice $Z$, and $z \in Z_+ \backslash A$. Then there exists $z^* \in Z^*_+$ so that $\langle z^*, z \rangle > \sup z^*|_A$.
\end{lem}

\begin{proof}
 By Hahn-Banach Theorem, there exists $x^* \in Z^*$ so that $\langle x^*, z \rangle > \sup x^*|_A$. For $z^* = x^*_+$, we have $\langle z^* , z \rangle \geq \langle x^* , z \rangle$. Further, by \cite[Lemma 19.8]{OiTu},
 $$
 \sup x^*|_A \geq \sup x^*|_{A \cap Z_+} = \sup z^*|_{A \cap Z_+} = \sup z^*|_A .
 $$
 Thus, $z^*$ has the desired properties.
\end{proof}

\subsection{The unit ball of $\fbp[E]$}\label{ss:fbp hull}

Begin by proving:

\begin{prop}\label{hull}
 If $E$ is a Banach space, and $1 \leq p \leq \infty$, then $\ball(\fbp[E])$ is the norm closed solid convex hull of 
 $$S = \Big\{ \big( \sum_{i=1}^N \big|\delta_{e_i}\big|^p \big)^{1/p} : N \in \Nat, \, e_1, \ldots, e_N \in E, \, \sum_i \|e_i\|^p \leq 1  \Big\}$$
 (when $p=\infty$, we use $S = \big\{ \vee_{i=1}^N \big|\delta_{e_i}\big| : N \in \Nat, \, e_1, \ldots, e_N \in E, \, \max_i \|e_i\| \leq 1  \big\}$ instead).
\end{prop}

\begin{proof}
 As $\fbp[E]$ is $p$-convex, the inequality
 $$ \big\| \big( \sum_{i=1}^N \big|\delta_{e_i}\big|^p \big)^{1/p} \big\| \leq  \big( \sum_i \big\|\delta_{e_i}\big\|^p \big)^{1/p} = \big( \sum_i \|e_i\|^p \big)^{1/p} $$
 holds for any $e_1, \ldots, e_N$. Therefore, set $S$, and consequently the closure of its convex hull, belong to $\ball({\fbp[E]})$.
 
Now suppose, for the sake of contradiction, that the inclusion described above is proper. By the ``solidity'' of the objects discussed here, a certain $f \in \ball(\fbp[E])_+$ lies outside of ${\mathrm{CSCH}}(S)$.  % $f \in \ball(\fbp[E])_+ \backslash S'$, where $S'$ is the closure of the solid convex hull of $S$. 
% ``Positive separation'' (\cite[Proposition 19.7]{OiTu}), combined with \Cref{separate}, allow us to
Combining Lemmas \ref{separate} and \ref{positive separation}, we find $\varphi \in \ball(\fbp[E])_+^*$ so that $\langle \varphi , f \rangle > 1 > \sup_{s \in S} \langle \varphi, s \rangle$. Equip $\fbp[E]$ with the semi-norm $\triple{x} = \langle \varphi, |x| \rangle$. Modding out by $\ker \triple{ \cdot }$ (which is an ideal), we obtain an $L$-space, which can be identified with $L_1(\mu)$, for some measure $\mu$. We shall denote by $J$ the natural lattice quotient from $\fbp[E]$ to $L_1(\mu)$. For future use, note that, if $e_1, \ldots, e_N \in E$ are such that $\sum_i \|e_i\|^p \leq 1$, then
\begin{equation}
\Big\| \big( \sum_i \big| J \delta_{e_i} \big|^p \big)^{1/p} \Big\|_{L_1(\mu)} = \big\langle \varphi , \big( \sum_i \big| \delta_{e_i} \big|^p \big)^{1/p} \big \rangle \leq 1 .
    \label{eq:J is nice}
\end{equation}

First assume that $p=1$. For any $e \in E$, we have $$\|J \circ \phi e\| = \|J \delta_e\| = \triple{\delta_e} \leq \|\delta_e\| = \|e\| , $$
hence $J\phi : E \to L_1(\mu)$ is a contraction (recall that $\phi : E \to \fbp[E]$ denotes the canonical embedding, and $\phi e$ is traditionally denoted by $\delta_e$). Therefore, $J\phi$ has a unique extension to a lattice homomorphism $\widetilde{J\phi} : \fbl[E] \to L_1(\mu)$. Then $\widetilde{J\phi} = J$, hence $\|\widetilde{J\phi}\| = \|J\| \leq 1$, leading to a contradiction: $1 < \|Jf\| \leq \|J\| \|f\| \leq 1$.

Next consider the case of $p \in (1,\infty)$. \eqref{eq:J is nice} guarantees that Maurey-Nikishin Factorization Theorem (see e.g.~\cite[Th{\'e}or{\`e}me 8]{Mau} or \cite[Proposition III.H.10]{Wojt}) is applicable to $J\phi$: there exist a contraction $u : E \to L_p(\nu)$, and a contractive lattice homomorphism $v : L_p(\nu) \to L_1(\mu)$ (in fact, implemented by a change of density), so that $J \phi = vu$. 
Then $u$ extends to a lattice homomorphism $\widetilde{u} : \fbp[E] \to L_p(\nu)$; $v \widetilde{u} : \fbp[E] \to L_1(\mu)$ is a lattice homomorphism as well, hence it coincides with $J$. Again we obtain a contradiction:
$1 < \|Jf\| \leq \|v\| \|\widetilde{u}\| \|f\| \leq 1$.

Finally, we come to $p = \infty$. The inequality $\| J \vee_{i=1}^N |\delta_{e_i}| \| = \|\vee_{i=1}^N J |\delta_{e_i}|\| \leq 1$ holds for any $e_1, \ldots, e_N \in \ball(E)$, hence, by the order completeness of $L_1(\mu)$, there exists a norm one $h \in L_1(\mu)_+$ so that $J |\delta_e| = |J \phi e| \leq h$ for any $e \in \ball(E)$. 
Consequently, we have a factorization $J\phi = vu$, where $v : L_\infty(\mu) \mapsto L_1(\mu) : g \mapsto hg$ is a lattice homomorphism, and $u : E \to L_\infty(\mu)$ is defined by $u e = h^{-1} \cdot J\phi e$ (if $h(t) = 0$, then $[J \delta_e](t) = 0$ as well; denote $0^{-1} \cdot 0 = 0$). Both $u$ and $v$ are contractive; now reason as for $p \in (1,\infty)$.
\end{proof}

We also recall the discussion, in \cite[Section 2]{OiTTT1}, on the \emph{free vector lattice} on $E$ (denoted by $\fvl[E]$). This is the (dense, by definition) sublattice of $\fbp[E]$, generated by $\phi(E)$ -- that is, the sublattice consisting of elements $F(\delta_{e_1}, \ldots, \delta_{e_N})$, with $e_1, \ldots, e_N \in E$, and $F$ a lattice expression (a formula using both linear and lattice operations; see e.g.~\cite[Section 1.d]{LT2} for positively homogeneous functional calculus). In fact, any element of $\fvl[E]$ can be written as $\vee_{i=1}^m \delta_{x_i} - \vee_{j=1}^n \delta_{y_j}$, for some $(x_i)$ and $(y_j)$.

The omission of any reference to the ambient $\fbp[E]$ in the notation $\fvl$ is justified: the proof of Theorem 2.1 from \cite{OiTTT1} shows that $F(\delta_{e_1}, \ldots, \delta_{e_N}) = G(\delta_{f_1}, \ldots, \delta_{f_M})$ iff $F(\langle e^*, e_1 \rangle, \ldots, \langle e^*, e_N \rangle) = G(\langle e^*, f_1 \rangle, \ldots, \langle e^*, f_M \rangle)$, for any $e^* \in E^*$; in this case, we say that $F(\delta_{e_1}, \ldots, \delta_{e_N})$ is \emph{trivially identical} to $G(\delta_{f_1}, \ldots, \delta_{f_M})$. Thus, for different values of $p$, the lattices generated in $\fbp[E]$ by $\phi_{E,p}(E)$ can be identified.

For future use, note that the same notion of a free vector lattice applies to free Banach lattices with upper $p$-estimates as well. Indeed, let $u$ and $v$ both denote the identity on $E$. These operators have unique extensions to contractive lattice homomorphisms $\overline{u} : \fbl[E] \to \fbl^{\uparrow p}[E]$ and $\overline{v} :\fbl^{\uparrow p}[E] \to \fbl^{(\infty)}[E]$. Suppose $F$ is a lattice expression, and $e_1, \ldots, e_N \in E$. Then, by \cite[Section 2]{JLTTT},
$$
\overline{u} F\big(\phi_{E,1}(e_1), \ldots, \phi_{E,1}(e_N)\big) = F\big(\phi_{E,\uparrow p}(e_1), \ldots, \phi_{E,\uparrow p}(e_N)\big) ,
$$
and
$$
\overline{v} F\big(\phi_{E,\uparrow p}(e_1), \ldots, \phi_{E,\uparrow p}(e_N)\big) = F\big(\phi_{E,\infty}(e_1), \ldots, \phi_{E,\infty}(e_N)\big) .
$$
Combined, these two identities show that $F(\delta_{e_1}, \ldots, \delta_{e_N}) = G(\delta_{f_1}, \ldots, \delta_{f_M})$ in $\fbl^{\uparrow p}[E]$ iff the two are trivially identical. Thus, we are entitled to talk about $\fvl[E]$ as a sublattice (dense, by \cite[Section 3]{JLTTT}) of $\fbl^{\uparrow p}[E]$.

\begin{cor}\label{c:one is enough}
 If $E$ is a Banach space, and $1 \leq p \leq \infty$, then, for any $\varphi \in \ball(\fbp[E])$, and any $\vr > 0$, there exist $e_1, \ldots, e_N \in E$ with $\sum_i \|e_i\|^p \leq 1$, and $\psi \in \fvl[E]$ so that $|\psi| \leq \big( \sum_{i=1}^N \big|\delta_{e_i}\big|^p \big)^{1/p}$, and $\|\varphi - \psi\| < \vr$. When $p=\infty$, we resort to the usual convention -- for instance, $\big( \sum_{i=1}^N \big|\delta_{e_i}\big|^p \big)^{1/p} = \vee_i \big|\delta_{e_i}\big|$.
\end{cor}

\begin{proof}
% For $p=\infty$, this statement is a straightforward coincidence of \Cref{hull}. We therefore focus on $p \in [1,\infty)$.
We begin by showing that, for any $\vr > 0$, there exist $e_1, \ldots, e_N \in E$ with $\sum_i \|e_i\|^p \leq 1$, and % $\chi \in \fbp[E]$ so that $|\chi| \leq \big( \sum_{i=1}^N \big|\delta_{e_i}\big|^p \big)^{1/p}$, and $\|\varphi - \chi\| < \vr$.
\begin{equation}
 \chi \in \fbp[E] \, \, \textrm{s.t.} \, \, |\chi| \leq \big( \sum_{i=1}^N \big|\delta_{e_i}\big|^p \big)^{1/p} \, \, \textrm{and} \, \, \|\varphi - \chi\| < \vr .
 \label{cond on chi}
\end{equation}

We present the proof for $p \in [1,\infty)$; the case $p = \infty$ is obtained by replacing the $p$-th roots of $p$-th powers by suprema.

 By \Cref{hull}, we can find $M \in \Nat$, $t_1, \ldots, t_M \geq 0$ with $\sum_j t_j = 1$, as well as collections $(e_{ij})_{i=1}^{N_j}$ ($1 \leq j \leq M$) with $\sum_i \|e_{ij}\|^p \leq 1$, as well as $\chi \in \fbp[E]$, so that $\|\chi - \varphi\| < \varepsilon$, and, for every $e^* \in E^*$, 
 $$|\chi(e^*)| \leq \sum_j t_j \big( \sum_i | \langle e^*, e_{ij} \rangle |^p \big)^{1/p}.$$
 Let $a_j(e^*) = \big( \sum_i | \langle e^*, e_{ij} \rangle |^p \big)^{1/p}$. By Jensen's Inequality, 
 $$
 |\chi(e^*)| \leq \sum_j t_j a_j(e^*) \leq \big( \sum_j t_j a_i(e^*)^p \big)^{1/p} = \Big( \sum_j t_j \sum_i | \langle e^*, e_{ij} \rangle |^p \Big)^{1/p} .
 $$
 Thus,
 $$
 |\chi| \leq \Big( \sum_{j,i} \big| \delta_{t_j^{1/p} e_{ij}} \big|^p \Big)^{1/p},
 \, {\textrm{  and  }} \,
 \sum_{j,i} \big\|t_j^{1/p} e_{ij}\|^p = \sum_j t_j \sum_i \|e_{ij}\|^p \leq 1 . % \qedhere
 $$
 
 Relabeling the vectors $(t_j^{1/p} e_{ij})$ into $(e_i)$, observe that we have find $\chi$ as in \eqref{cond on chi}.
  It remains to replace $\chi$ by $\psi \in \fvl[E]$. Fix $\sigma > 0$.  By the density of $\fvl[E]$, we can find disjoint $\eta, \xi \in \fvl[E]_+$ so that $\|\chi_+ - \eta\|  , \|\chi_- - \xi\| < \sigma$ (indeed, begin by approximating $\chi_+, \chi_-$ by $\eta_0, \xi_0$, and then pass to $\xi = \xi_0 - \xi_0 \wedge \eta_0$, $\eta = \eta_0 - \xi_0 \wedge \eta_0$). 
  
  As $\ball(\ell_p^N)$ can be approximated with symmetric polytopes arbitrarily well, we can find a lattice expression (that is, a combination of linear and lattice operations) $\Phi : \Real^N \to \Real$ so that, for any $t_1, \ldots, t_N \in \Real$, we have
  $$ \Phi(t_1, \ldots, t_N) \leq \Big( \sum_i |t_i|^p \Big)^{1/p} \leq (1+\sigma) \Phi(t_1, \ldots, t_N) . $$
  Let $\zeta = \Phi \big ( \big| \delta_{e_1} \big| , \ldots ,  \big| \delta_{e_N} \big| )$.
  Then $\zeta \in \fvl[E]$,   
  $$\zeta \leq \big( \sum_{i=1}^N \big|\delta_{e_i}\big|^p \big)^{1/p} \leq (1+\vr) \zeta , \, \, {\textrm{hence}} \, \, \big\| \big( \sum_{i=1}^N \big|\delta_{e_i}\big|^p \big)^{1/p} - \zeta \big\| \leq \sigma . $$
  
  Let $\psi_+ = \eta \wedge \zeta$, $\psi_- = \xi \wedge \zeta$ (these two are disjoint), and $\psi = \psi_+ - \psi_-$. Clearly $\psi = \psi_+ - \psi_-$ lies in $\fvl[E]$, and $|\psi| \leq \zeta \leq \big( \sum_{i=1}^N \big|\delta_{e_i}\big|^p \big)^{1/p}$. To show that $\psi$ provides a good approximation for $\varphi$, observe first that, from the definition of $\psi_+$,
  $$
  |\eta - \psi_+| = (\eta - \zeta)_+ \leq \sigma \zeta + \Big( \eta - \eta \wedge (1+\sigma)\zeta \Big)_+ \leq \sigma \zeta + \Big( \eta - \eta \wedge \Big( \sum_i \big|\delta_{e_i}\big|^p \Big)^{1/p} \Big)_+ ,
  $$
  and further, as $\chi_+ \leq \big( \sum_i \big|\delta_{e_i}\big|^p \big)^{1/p}$,
  $$
  \Big( \eta - \eta \wedge \Big( \sum_i \big|\delta_{e_i}\big|^p \Big)^{1/p} \Big)_+ \leq |\eta - \chi_+| .
  $$
  Thus,
  $$
  \|\psi_+ - \chi_+\| \leq \|\chi_+ - \eta\| + \|\eta - \psi_+\| \leq 2 \|\chi_+ - \eta\| + \sigma \|\zeta\| \leq 3 \sigma ,
  $$
and likewise, $\|\psi_- - \chi_-\| \leq 3 \sigma$. Thus, $\|\psi - \chi\| \leq 6 \sigma$. By selecting $\sigma$ to be sufficiently small, we guarantee $\|\varphi - \psi\| \leq \|\varphi - \chi\| + \|\psi - \chi\| < \vr$.
\end{proof}

% From this we conclude:

The preceding result helps us to uncover a link between elements of $\fbp[E]$ and $p$-summing operators.
Recall that $T \in B(E,F)$ ($E,F$ are Banach spaces) is \emph{$p$-summing} if there exists a constant $C$ so that $\big( \sum_i \|T x_i\|^p \big)^{1/p} \leq C \| (x_i) \|_{p, {\textrm{weak}}}$, for any $x_1, \ldots, x_N \in E$.
The infimum of all $C$'s as above is called the \emph{$p$-summing norm} of $T$, and denoted by $\pi_p(T)$.
The famous Pietsch Factorization Theorem states that $\pi_p(T) \leq C$ if an only if it has a factorization $T = V \circ J|_{\overline{\ran U}} \circ U$, with $\|U\| \|V\| \leq C$; here $U \in B(E, L_\infty(\mu))$ ($\mu$ is a probability measure), $J : L_\infty(\mu) \to L_p(\mu)$ is the formal identity, and $V \in B(\overline{\ran JU},F)$.
For more information about $p$-summing operators, we refer the reader to e.g.~\cite{DJT}.

\begin{prop}\label{p:fvl in fbp}
 If $E$ is a Banach space, and $1 \leq p < \infty$, then, for $\varphi \in \fvl[E]$, the following are equivalent:
 \begin{enumerate}
  \item\label{norm one} $\|\varphi\|_p \leq 1$. % (here $\| \cdot \|_p$ stands for the norm in $\fbp[E]$).
  \item\label{inf many} For any $\varepsilon > 0$, there exist $e_1, e_2, \ldots \in E$ so that $\sum_i \|e_i\|^p < 1 + \varepsilon$, and the inequality $|\varphi(e^*)|^p \leq \sum_i | \langle e^*, e_i \rangle|^p$ holds for any $e^* \in E^*$.
  \item\label{p-summing} There exist a Banach space $Z$, and an operator $T : E^* \to Z$ with $\pi_p(T) \leq 1$, so that $|\varphi(e^*)| \leq \|T e^*\|$ for any $e^* \in E^*$; in fact, we can take $Z = L_p(\nu)$, for some measure $\nu$.
  \item\label{dominant measure} There exists a Radon probability measure $\mu$ on $\ball(E^{**})$ so that $|\varphi(e^*)|^p \leq \int | \langle e^*, e^{**} \rangle|^p \, d\mu(e^{**})$ for any $e^* \in E^*$.
 \end{enumerate}
\end{prop}

\begin{proof}
The equivalence \eqref{norm one} $\Leftrightarrow$ \eqref{dominant measure} has appeared in \cite{JLTTT}; we prove all implications for the reader's convenience.

\eqref{norm one} $\Rightarrow$ \eqref{inf many}: without loss of generality, we can assume $\varphi \geq 0$. Let $t_0 = 1 = s_0$, and, for $k \geq 1$, fix positive $t_k, s_k$ so that
$$
\sum_{k=0}^\infty t_k^{p'} < (1+\vr)^{p'/2} , \, \, \sum_{k=0}^\infty s_k^p < (1+\vr)^{p/2}. % , \, \, {\textrm{with}} \, \, \frac1p + \frac1{p'} = 1.
$$
By \Cref{c:one is enough}, we can find $\varphi_0, \varphi_1, \ldots \in \fvl[E]_+$ so that:
\begin{enumerate}
 \item $\|\varphi_k\| < t_k s_k$ for $k \geq 0$.
 \item For any $k \geq 0$, $\|\varphi - \sum_{i=0}^k \varphi_i\| < t_{k+1} s_{k+1}$.
 \item For each $k$, there exists $e_{k,1}, \ldots, e_{k,M_k} \in E$, so that $\sum_i \|e_{ki}\|^p \leq s_k^p$, and $\varphi_k \leq t_k \psi_k$, where $\psi_k = \big( \sum_i \big| \delta_{e_{ki}} \big|^p \big)^{1/p}$.
\end{enumerate}
Note that
$$
\sum_k \sum_i \|e_{k,i}\|^p < \sum_k s_k^p = (1 + \vr)^{p/2} 
$$
Then, by H\"older Inequality,
$$
\varphi \leq \sum_{k=0}^\infty t_k \psi_k \leq \big( \sum_k t_k^{p'} \big)^{1/p'} \big(\sum_k \psi_k^p \big)^{1/p} < \sqrt{1+\vr} \big(\sum_k \psi_k^p \big)^{1/p} .
$$
so, for any $e^* \in E^*$, we have
$$
\varphi(e^*) \leq \sqrt{1+\vr} \big( \sum_{k,i} | \langle e^*, e_{k,i} \rangle |^p \big)^{1/p} .
$$
It remains to recall that
$$
\sum_k \sum_i \|e_{k,i}\|^p < \sum_k s_k^p = (1 + \vr)^{p/2} .
$$

\eqref{inf many} $\Rightarrow$ \eqref{p-summing}: For each $N \in \Nat$, \eqref{inf many} permits us to find non-zero $(e_{i,N})_{i \in \Nat} \subset E$ so that $\sum_i \|e_{i,N}\|^p < (1+2^{-N})^p$, and $|\varphi(e^*)|^p \leq \sum_i |\langle e^*, e_{i,N} \rangle |^p$ for any $e^* \in E^*$.
Then the operator $T_N : E^* \to \ell_p : e^* \mapsto (\langle e^*, e_{i,N} \rangle)_i$ can be written as $T_N = \Delta S$, where the contraction $S : E^* \to \ell_\infty$ takes $e^*$ to $ (\|e_{i,N}\|^{-1} \langle e^*, e_{i,N} \rangle)_i$, and $\Delta : \ell_\infty \to \ell_p$ is a diagonal operator, multiplying the $i$-th diagonal entry by $\|e_{i,N}\|$. So, $\pi_p(T_N) < 1+2^{-N}$.

Now suppose ${\mathcal{U}}$ is a free ultrafilter on $\Nat$. Then $\widetilde{T} := \prod_{\mathcal{U}} T_N$ has some $L_p(\mu)$ as its range, and satisfies $\pi_p(\widetilde{T}) \leq 1$. Then $T := \widetilde{T}\big|_{E^*}$ has the desired properties.

 \eqref{p-summing} $\Rightarrow$  \eqref{dominant measure} is a direct consequence of Pietsch Factorization Theorem.

 \eqref{dominant measure} $\Rightarrow$ \eqref{norm one}: If $e_1^*, \ldots, e_n^* \in E^*$ satisfy $\| (e_j^*) \|_{p, {\textrm{weak}}} \leq 1$, then
 $$
 \sum_j |\varphi(e_j^*)|^p \leq \sum_j \int | \langle e_j^*, e^{**} \rangle|^p \, d\mu(e) = \int \sum_j | \langle e_j^*, e^{**} \rangle|^p \, d\mu(e^{**}) \leq \| (e_j^*) \|_{p, {\textrm{weak}}}^p \leq 1 ,
 $$
 yielding $\|\varphi\|_p \leq 1$.
\end{proof}

\subsection{The unit ball of $\fbl^{\uparrow p}[E]$}\label{ss:fb upper p hull}

Below we present an analogue of \Cref{hull} for $\fbl^{\uparrow p}[E]$.

\begin{prop}\label{hull upper est}
 Fix $p \in (1,\infty)$, and a Banach space $E$. Let $S$ be the subset of $\fbl^{\uparrow p}[E]$, consisting of functions $\vee_i |\delta_{e_i}|$, with $\sum_i \|e_i\|^p \leq 1$; let $S'$ be the closure of the solid convex hull of $S$. Then:
 \\
 $(1)$ $S' \subset \ball(\fbl^{\uparrow p}[E])$.
 \\
 $(2)$ There exists $K = K_p$ so that $\ball(\fbl^{\uparrow p}[E]) \subset K_p S'$.
\end{prop}

% Before presenting the proof, we make some remarks.

% No good functional representation for $\fbl^{\uparrow p}[E]$ seems to be known. It is easy to see that there exists a 1-1 isometric correspondence between elements of $E^*$ and lattice homomorphisms $\fbl^{\uparrow p}[E] \to \Real$. However, it is not clear (to me) that such lattice homomorphisms separate points; thus we may not be able to view $\fbl^{\uparrow p}[E]$ as a set of homogeneous functions on $E^*$.

% If $E$ is finite dimensional, then we can identify identify $\fbl^{\uparrow p}[E]$ with $\sphere(E^*)$, equipped with the norm described by \Cref{hull upper est} (the norm in question is equivalent to $\| \cdot \|_\infty$, so the closure can be taken with respect to that norm as well.

The proof given below resembles that of \Cref{hull}, but, instead of Maurey-Nikishin Theorem for factoriations via $L_p$, we need to rely on Pisier's result on operators factoring via weak $L_p$. 
 Specifically, by the proof of \cite[Theorem 1.2]{Pisier_L_p_infty}, there exists a constant $K_p$ with the following property: suppose $Z$ is a Banach space, and $T : Z \to L_1(\mu)$ is an operator with the property that, for any $z_1, \ldots, z_n \in Z$, we have $\| \vee_i |Tz_i| \| \leq (\sum_i \|z_i\|^p)^{1/p}$. Then $T$ can be written as $vu$, where $\nu$ is a probability measure, $v : L_{p\infty}(\nu) \to L_1(\mu)$ is a contractive lattice homomorphism (a change of density), and $u \in B(Z, L_{p\infty}(\nu))$ satisfies $\|u\| \leq K_p$.
 Here we equip $L_{p\infty}(\Omega,\nu)$ with the norm
 \begin{equation}
 \|h\| = \|h\|_{p,\infty} = \sup_{E \subset \Omega} \nu(E)^{1/p - 1} \int_E |h| .
 \label{eq:p infty new norm}
 \end{equation}
 This norm turns $L_{p\infty}(\Omega,\nu)$ into a Banach lattice, which satisfies the upper $p$-estimate with constant $1$ \cite[Section 9.6]{OiTTT1}.
 By e.g.~\cite[Exercise 1.1.12]{Grafakos}, if $\nu$ is $\sigma$-finite, then $\| \cdot \|_{p,\infty}$ is equivalent to the ``canonical'' measure on $L_{p\infty}$:
\begin{equation}
 \triple{h}_{p,\infty} = \sup_{t > 0} t^{1/p} \nu(|h| > t) 
\label{eq:p infty norm} \end{equation}
(the constant of equivalence depends on $p$ only).
 
\begin{proof}[Proof of \Cref{hull upper est}]
(1) As $\fbl^{\uparrow p}[E]$ has upper $p$-estimate with constant $1$, we have $\|\vee_i |\delta_{e_i}| \| \leq 1$ whenever $\sum_i \|e_i\|^p \leq 1$.
 
 (2) Suppose, for the sake of contradiction, that there exists $f \in \ball(\fbl^{\uparrow p}[E]) \backslash K_p S'$. By solidity, we can assume $f \geq 0$. Find $\varphi \in \fbl^{\uparrow p}[E]^*_+$ which separates $f$ from $S'$ -- that is, $\langle \varphi,f \rangle >K_p$ and $1 > \sup_{s \in S'} \langle \varphi , s \rangle$. As in the proof of \Cref{hull}, let, for $x \in \fbl^{\uparrow p}[E]$, $\triple{x} = \langle \varphi, |x| \rangle$, and identify $(\fbl^{\uparrow p}[E] , \triple{ \cdot }) /\ker \triple{ \cdot }$ with an $L_1(\mu)$, for some measure $\mu$.
 Denote the canonical embedding of $\fbl^{\uparrow p}[E]$ into $L_1(\mu)$ by $J$ (this is clearly a lattice homomorphism). Let $\phi$ be the canonical embedding of $E$ into $\fbl^{\uparrow p}[E]$. If $\sum_{i=1}^n \|e_i\|^p \leq 1$, then
$$
\big\| \vee_i |J \phi e_i| \big\| = \big\| \vee_i |\delta_{e_i}| \big\| = \triple{ \vee_i |J \delta_{e_i}|} \leq 1 ,
$$
and therefore, as explained above, we can factor $J \phi : E \to L_1(\mu)$ as $vu$, were $v : L_{p\infty}(\nu) \to L_1(\mu)$ is a contractive lattice homomorphism, and $\|u\| \leq K_p$. Let $\widehat{u} : \fbl^{\uparrow p}[E] \to L_{p\infty}(\nu)$ be the lattice homomorphism canonically extending $u$; then $v \widehat{u} = \widehat{J \phi}$ has norm not exceeding $K_p$, which contradicts the assumption $\triple{f} = \|Jf\| > K_p$.
\end{proof}

\section{Positive Bounded Approximation Property}\label{bap}

Begin by introducing a definition. 
% Suppose $E$ is a finite dimensional Banach space. We say that $\| \cdot \|$ is a \emph{convenient} norm on $C(\sphere(E^*))$
For a finite dimensional Banach space $E$, the norm $\| \cdot \|$ on $C(\sphere(E^*))$ is called \emph{convenient} if the following conditions are satisfied:
\begin{enumerate}
 \item $\| \cdot \|$ is a lattice norm (with pointwise operations).
 \item There exists $c > 0$ s.t.~$c \| \cdot \| \leq \| \cdot \|_\infty \leq \| \cdot \|$.
 \item There exists an equicontinuous family $\css \subset C(\sphere(E^*))_+$ (the \emph{controlling family}) so that $\|f\| \leq 1$ if and only if for every $\vr > 0$ there exists $g \in \css$ so that $\| |f| \vee g - g\|_\infty \leq \vr$.
\end{enumerate}
Recall that the equicontinuity of $\css$ means that there exists a non-decreasing function $\boldom = \boldom_{\css} : (0,\infty) \to (0,\infty)$, so that $\lim_{s \to 0+} \boldom(s) = 0$, and $|g(x_1^*) - g(x_2^*)| \leq \boldom(\|x_1^* - x_2^*\|)$ for any $x_1^*, x_2^* \in \sphere(E^*)$.
By Ascoli-Arzela Theorem, the closure of such $\css$ is compact.
% for every $\vr > 0$ there exists $\delta > 0$ so that any $g \in \css$ satisfies $|g(x_1^*) - g(x_2^*)| < \vr$ whenever $\|x_1^* - x_2^*\| < \delta$.

\begin{rem}\label{r:examples of convenience}
As before, suppose $E$ is finite dimensional.
%  (1) The prime example of a convenient norm is that of $\fbl[E]$. In this case, by \cite{OiTT2}, $(\dim E)^{-1} \| \cdot \|_\infty \leq \| \cdot \| \leq \| \cdot \|_\infty$. In this case, by Theorem \ref{t:hull of extremes}, the controlling family $\css$ consists of all finite sums $g = \sum_i |\delta_{e_i}|$, with $\sum_i \|e_i\| = 1$. Functions $g$ as above are $1$-Lipschitz. Indeed, for $x_1^*, x_2^* \in \sphere(E^*)$,
%  \begin{align*}
%  |g(x_1^*) - g(x_2^*)|   & =  \big| \sum_i \big( | \langle x_1^*, e_i \rangle | - | \langle x_2^*, e_i \rangle | \big) \big| 
%  \\ & \leq \sum_i \big| \langle x_1^*, e_i \rangle - \langle x_2^*, e_i \rangle  \big| \leq \sum_i \|x_1^* - x_2^*\| \|e_i\| = \|x_1^* - x_2^*\| ,
%  \end{align*}
%  hence $\css$ is equicontinuous.
 %
%  (2) More generally, the norm of $\pfbp[E]$ ($1 \leq p < \infty$) is convenient. By Corollary \ref{c:pfbp hull}, we can take as $\css$ the set of finite sums $g = \sum_i \lambda_i |\delta_{e_i}|^p$, with $\lambda_i \geq 0$, $\sum_i \lambda_i = 1$, and $e_i \in \sphere_E$. To establish equicontinuity, consider $x_1^*, x_2^* \in \sphere_{E^*}$, and $e \in \sphere_E$. Mean Value Theorem shows that, for $t, s \in [0,1]$, we have $|t^p - s^p| \leq p |t-s|$, and therefore,
%  $$  |\langle x_1^*, e\rangle|^p - |\langle x_2^*, e\rangle|^p \leq p \big| |\langle x_1^*, e\rangle| - |\langle x_2^*, e\rangle| \big| \leq p \|x_1^* - x_2^*\| . $$
%  Therefore, $|g(x_1^*) - g(x_2^*)| \leq p \|x_1^* - x_2^*\|$ for any $g \in \css$, and $x_1^*, x_2^* \in \sphere_{E^*}$.
 %
 
 (i) The norm of $\fbp[E]$ is convenient. Indeed, it is a lattice norm, so (1) holds. To show (2), note that, for any $f \in C(\sphere(E^*))$, we have $\|f\| \geq \|f\|_\infty$ (here $\| \cdot \|$ is the $\fbp[E]$-norm). 
 To estimate $\|f\|$ from above, fix $x_1^*, \ldots, x_N^* \in E^*$ with $\|(x_k^*)\|_{p,{\textrm{weak}}} \leq 1$. Then $\sum_k f(x_k^*)^p \leq \|f\|_\infty^p \sum_k \|x_k^*\|^p \leq \pi_p(id_E)^p \|f\|_\infty^p$. However, it is well know (see e.g.~\cite{TJ}) that $\pi_p(id_E) \leq \dim E$. Recalling the description of $\fbp[E]$ given in the end of \Cref{intro}, we conclude that $\|f\| \leq \dim E \cdot \|f\|_\infty$.
 
 Finally we verify (3). For $p = \infty$, it suffices to take $\css = \{ 1 \}$ (constant function equal $1$ everywhere). For $1 \leq p < \infty$, use $\css$ consisting of all $g = ( \sum_{i=1}^N \lambda_i |\delta_{e_i}|^p )^{1/p}$ with $\lambda_1, \ldots, \lambda_N \in [0,1]$, $\sum_i \lambda_i = 1$, and $e_1, \ldots, e_N \in \sphere(E)$. 
 The mean value theorem shows that, for $t, s \in [0,1]$, we have $|t^p - s^p| \leq p |t-s|$. Therefore, for $x_1^*, x_2^* \in \sphere(E^*)$,
 $$
 \big| |\langle x_1^*, e_i \rangle|^p - |\langle x_2^*, e_i \rangle|^p \big| \leq p \big| |\langle x_1^*, e_i \rangle| - |\langle x_1^*, e_i \rangle| \big| \leq p |\langle x_1^* - x_2^*, e_i \rangle| \leq p \|x_1^* - x_2^*\| ,
 $$
 and therefore,
 $$
 \big| g(x_1^*)^p - g(x_2^*)^p \big| \leq \sum_i \lambda_i \big| |\langle x_1^*, e_i \rangle|^p - |\langle x_2^*, e_i \rangle|^p \big| \leq p \|x_1^* - x_2^*\| .
 $$
%  Note also that $\|g\|_\infty \leq \|g\| \leq 1$, hence $|g(x_1^*)|, |g(x_2^*)| \leq 1$.
 Recall that $|a-b|^p \leq |a^p - b^p|$ for $a,b \geq 0$; therefore,
 $$
 |g(x_1^*) - g(x_2^*)| \leq \big| g(x_1^*)^p - g(x_2^*)^p \big|^{1/p} \leq p^{1/p} \|x_1^* - x_2^*\|^{1/p} ,
 $$
 establishing the uniform continuinty of $\css$.
 
 It remains to show that $\css$ determines the norm of $\fbp[E]$ -- that is, $f \in C(\sphere(E^*))_+$ satisfies $\|f\| \leq 1$ iff for every $\vr > 0$ there exists $g \in \css$ so that $\|f \vee g - f\|_\infty < \vr$. By \Cref{c:one is enough}, there exists $g \in \css$ and $h \in \fbp[E]$ so that $0 \leq h \leq g$, and $\|f-h\| \leq \vr$. Then
 $g + (f-h)_+ \geq h + (f-h) = f$, hence $f \vee g \leq g + (f-h)_+$, and consequently,
 $\|f \vee g - f\|_\infty \leq \|f-h\| < \vr$.
 
 (ii) For $1 < p < \infty$, the formal identities $\fbl[E] \to \fbl^{\uparrow p}[E] \to \fbl^{(\infty)}[E]$ are contractive, hence, in light of part (i) we can view $\fbl^{\uparrow p}[E]$ as a space of continuous functions on $\sphere(E^*)$. % with the norm equivalent to $\| \cdot \|$.
 Below we show that there exists a convenient norm $\triple{ \cdot }$, so that $\triple{ \cdot } \leq \| \cdot \| \leq K_p \triple{ \cdot }$ (here $\| \cdot \|$ is the usual norm of $\fbl^{\uparrow p}[E]$, and $K_p$ is the constant from \Cref{hull upper est}).
 
 Let $\css$ be the family consisting of convex combinations of functions $\vee_{i=1}^n |\delta_{e_i}|$, with $\sum_i \|e_i\|^p \leq 1$. This family is equicontinuous. Indeed, it suffices to show that, for any $g \in \css$ as above, and any $x_1^*, x_2^* \in \sphere(E^*)$, we have $g(x_2^*) \geq g(x_1^*) - \|x_1^* - x_2^*\|$. 

 Write $g = \sum_{j=1}^N t_j g_j$, where $t_j \geq 0$, $\sum_j t_j = 1$, and $g_j = \vee_{i=1}^{n_j} |\delta_{e_{ij}}|$.
 For each $j$, find $i$ so that $g_j(x_1^*) = |\langle x_1^*, e_{ij} \rangle|$. As $\|e_{ij}\| \leq 1$, 
 \begin{align*}
 g_j(x_2^*) 
 &
 \geq |\langle x_2^*, e_{ij} \rangle| = |\langle x_1^* + (x_2^* - x_1^*) , e_{ij} \rangle| \geq |\langle x_1^*, e_{ij} \rangle| - | \langle x_2^* - x_1^* , e_{ij} \rangle|
 \\
 &
 \geq g_j(x_1^*) - \|x_1^* - x_2^*\|.  
 \end{align*}
 Adding these inequalities with weights $t_j$, we obtain $g(x_2^*) \geq g(x_1^*) - \|x_1^* - x_2^*\|$.
 
 Let $\triple{ \cdot }$ be the norm generated by the equicontinuous family $\css$. By \Cref{hull upper est}, $\triple{ \cdot } \leq \| \cdot \| \leq K_p \triple{ \cdot }$. Above we have verified part (3) of the definition of convenience. Clearly $\triple{ \cdot }$ is a lattice norm, hence part (1) of that definition also holds. 
 As this norm is dominated by that of $\fbl[E]$, there exists a constant $c$ so that $c \triple{ \cdot } \leq \| \cdot \|_\infty$; this handles the left hand side of condition (2). 
 Finally, $\|g\|_\infty \leq 1$ for any $g \in \css$, hence the right side of condition (2) holds as well.
 
 We do not know whether the ``original'' norm $\| \cdot \|$ of $\fbl^{\uparrow p}[E]$ is convenient.
\end{rem}

\begin{thm}\label{t:general PMAP}
If $\| \cdot \|$ is a convenient norm on $C(\sphere(E^*))$, then $\big(C(\sphere(E^*)), \| \cdot \| \big)$ has the $1$-PAP.
\end{thm}

In light of Remark \ref{r:examples of convenience}, we immediately conclude:

\begin{cor}\label{c:examples of PMAP}
Suppose $E$ is finite dimensional.
\begin{enumerate}
 \item For any $p \in [1,\infty]$, $\fbp[E]$ has the $1$-PAP.
 \item For any $p \in (1,\infty)$, $\fbl^{\uparrow p}[E]$ has the $K_p$-PAP. 
\end{enumerate}
%  If $E$ is finite dimensional, then, for any $p \in [1,\infty]$, $\fbp[E]$ has the PMAP. % In particular, $\fbl[E]$ has the PMAP.
\end{cor}

% \begin{cor}\label{c:example of PBAP}
%  If $E$ is finite dimensional, then, for any $p \in (1,\infty)$, $\fbl^{\uparrow p}[E]$ has the $K_p$-PAP. % In particular, $\fbl[E]$ has the PMAP.
% \end{cor}

\begin{proof}[Proof of Theorem \ref{t:general PMAP}]
For brevity, denote $\big(C(\sphere(E^*)), \| \cdot \| \big)$ by $X$. Our goal is to find a net of positive finite rank operators $(u_\alpha)$ on $X$, so that $u_\alpha \to 1$ point-norm, and $\limsup_\alpha \|u_\alpha\| \leq 1$.
 
%  Let $\mathcal{A}$ be the set of all partitions of unity of $\sphere_{E^*}$ -- that is, of all representations $1_{\sphere_{E^*}} = \sum_i f_{i \alpha}$, where $f_{i \alpha}$ are non-negative functions, attaining the values $1$

Let ${\mathcal{A}}$ be the set of pointed partitions of unity -- that is, of all families $\alpha = (x_1^*, \ldots, x_n^*, f_1, \ldots, f_n)$, where:
\begin{enumerate}
 \item $x_1^*, \ldots, x_n^* \in \sphere(E^*)$, $f_1, \ldots, f_n \in C(\sphere(E^*))_+$;
 \item $\sum_i f_i = 1$, $f_i(x_i^*) = 1$ for $i = 1, \ldots, n$.
\end{enumerate}
Let $\diam \alpha = \max_i \diam \supp f_i$. We say $\alpha \prec \beta$ if $\diam \beta \leq \diam \alpha$.

For $\alpha$ as above, and $f \in C(\sphere(E^*))$, define $P_\alpha f = \sum_{i=1}^n f(x_i^*) f_i$. Clearly $P_\alpha \geq 0$. By uniform continuity, $P_\alpha f \to f$ in the $\| \cdot \|_\infty$ norm, hence also in $\| \cdot \|$ (since the two norms are equivalent). Also,  $\|P_\alpha\|_\infty \leq 1$.

Fix $\vr > 0$, and $f \in C(\sphere(E^*))_+$ with $\|f\| < 1$. We shall show that
\begin{equation}
\label{eq:P alfa f}
\|P_\alpha f\| \leq 1 + c^{-1} \boldom(\diam \alpha) + \vr .
\end{equation}
Recalling that $\vr$ can be arbitrarily close to $0$, and invoking the positivity of $P_\alpha$, we conclude that $\|P_\alpha\| \leq 1 + c^{-1} \boldom(\diam \alpha)$.

Find $g \in \css$ and $h \in C(\sphere(E^*))_+$ so that $\|h\|_\infty < c \vr$, and $f \leq g + h$. Then $P_\alpha f \leq P_\alpha g + P_\alpha h$, and $\|P_\alpha h\| \leq c^{-1} \|P_\alpha h\|_\infty < \vr$. Consequently, \eqref{eq:P alfa f} will be established once we show that 
\begin{equation}
\|P_\alpha g - g\|_\infty \leq \boldom(\diam \alpha) 
\label{eq:difference}
\end{equation}
(it would then follow that $\|P_\alpha g - g\| \leq c^{-1} \boldom(\diam \alpha)$).
% there exists $u \in C(\sphere_{E^*})_+$ so that $\|u\|_\infty \leq \boldom(\diam \alpha)$ (so $\|u\| \leq c^{-1} \boldom(\diam \alpha)$), and $P_\alpha g \leq g + u$.

To establish \eqref{eq:difference}, fix $x^* \in \sphere(E^*)$.
% and show that
% \begin{equation}
% |[P_\alpha g](x^*) - g(x^*)|_\infty \leq ??? .
% \label{eq:pointwise difference}
% \end{equation}
Let $S = \{i : f_i(x^*) > 0\}$, then $g(x^*) = \sum_{i \in S} f_i(x^*) g(x^*)$ (since $\sum_{i \in S} f_i(x^*) = 1$), and $[P_\alpha g](x^*) = \sum_{i \in S} f_i(x^*) g(x_i^*)$. As $\|x^* - x_i^*\| \leq \diam \alpha$ for $i \in \css$, we have $|g(x^*) - g(x_i^*)| \leq \boldom(\diam \alpha)$, and therefore,
$$
|[P_\alpha g](x^*) - g(x^*)| \leq \sum_{i \in S} f_i(x^*) |g(x_i^*) - g(x^*)| \leq \boldom(\diam \alpha) 
$$
(here we again use $\sum_{i \in S} f_i(x^*) = 1$). As $x^* \in \sphere(E^*)$ is arbitrary, \eqref{eq:difference} follows.
\end{proof}

\begin{thm}\label{t:bap to pbap} 
 If a Banach space $E$ has the $\lambda$-Approximation Property, then $\fbp[E]$ has $\lambda$-Positive Approximation Property.
\end{thm}

\begin{proof}
 Suppose $(u_\alpha)$ is a net of finite rank operators, converging to $I_E$ point-norm, and such that $\sup_\alpha \|u_\alpha\| \leq \lambda$. Let $E_\alpha = \ran(u_\alpha)$, and let $j_\alpha : E_\alpha \to E$ be the corresponding embedding. % We shall identify $u_\alpha$ with its astriction to $E_\alpha$.
 Abusing the notation slightly, we shall view $u_\alpha$ as an operator from $E$ to $E_\alpha$.
 
 ``Liberation'' produces lattice homomorphisms $\overline{u_\alpha} : \fbp[E] \to \fbp[E_\alpha]$ and $\overline{j_\alpha} : \fbp[E_\alpha] \to \fbp[E]$ (see \cite{OiTTT1} for the properties of the extension of $u$ to $\overline{u}$). Note that $j_\alpha u_\alpha \to I_E$ point-norm, hence, by the continuity of lattice operations, $\overline{j_\alpha u_\alpha} f \to f$ for any $f \in \fvl[E]$ (the free vector lattice on $E$, which we view as the set of functions represented as finite lattice expressions in elements on $\sphere(E^*)$). As $\overline{j_\alpha u_\alpha} = \overline{j_\alpha} \overline{u_\alpha}$ has norm not exceeding $\lambda$, and $\fvl[E]$ is dense in $\fbp[E]$ \cite{OiTTT1}, we conclude that $\overline{j_\alpha} \overline{u_\alpha} \to I_{\fbl[E]}$ point-norm.
 
 Now take a finite dimensional subspace $Z \subset \fbp[E]$, and $\varepsilon > 0$. For $\alpha$ large enough, $\|\overline{j_\alpha} \overline{u_\alpha} z - z\| < \varepsilon$ for any $z \in \ball(Z)$. Further, $\overline{u_\alpha} Z$ is finite dimensional, hence there exists a positive finite rank contraction $v \in B(\fbp[E_\alpha])$, so that $\|(v-I)|_{\overline{u_\alpha} Z}\| < \vr/\lambda$.
 Then $\overline{j_\alpha} \, v \, \overline{u_\alpha}$ is a finite rank positive operator of norm not exceeding $\lambda$. Moreover, for any $z \in \ball(Z)$, 
 \begin{align*}
 \|\overline{j_\alpha} v \overline{u_\alpha}z - z\| 
 &
 \leq \|\overline{j_\alpha} (v-I) \overline{u_\alpha}z\| + \|\overline{j_\alpha} \overline{u_\alpha} z - z\| 
 \\
 &
 < \|(v-I)|_{\overline{u_\alpha} Z}\| \|\overline{u_\alpha} z\| + \vr < (\lambda + 1) \vr ,
 \end{align*}
 hence $ \|\overline{j_\alpha} v \overline{u_\alpha}|_Z - I_Z\| < (\lambda+1)\varepsilon$.
\end{proof}

In a similar fashion we show:

\begin{thm}\label{t:bap for upper p est} 
 If a Banach space $E$ has the $\lambda$-Approximation Property, and $1 < p < \infty$, then $\fbl^{\uparrow p}[E]$ has $\lambda K_p$-Positive Approximation Property.
\end{thm}

It is an open problem whether any Banach lattice with the $\lambda$-AP must also have $\nu$-PAP, for some $\nu$.
To partially resolve this, we introduce a new definition: we say that a Banach lattice $X$ \emph{$\lambda$-positively factors} through $\fbp[X]$ if there exist positive maps
$u : X \to \fbp[X]$ and $v : \fbp[X] \to X$ so that $I_X = vu$, and $\|u\| \|v\| \leq \lambda$. 
$\lambda$-positive factorization through $\fbl^{\uparrow p}[E]$ is defined similarly.

Theorem \ref{t:bap to pbap} immediately implies:

\begin{cor}\label{c:bap to pbap general}
 If a Banach lattice $X$ has the $\mu$-Approximation Property, and $\lambda$-positively factors through $\fbp[X]$ ($\fbl^{\uparrow p}[X]$), then $X$ has the $\lambda\mu$-Positive Approximation Property (respectively, $K_p\lambda\mu$-Positive Approximation Property, with $K_p$ as in \Cref{c:examples of PMAP}).
\end{cor}

In \cite{God-Kal Principle}, a stronger property was considered: a Banach lattice $X$ is said to have the \emph{$\lambda$-lattice lifting property} (through $\fbp[X]$ or $\fbl^{\uparrow p}[X]$) if, in the definition of positive factorization, $u$ and $v$ can be taken to be lattice homomorphisms. By \cite[Proposition 2.3]{God-Kal Principle}, for any Banach space $E$, $\fbp[E]$ ($\fbl^{\uparrow p}[E]$) has the $1$-lattice lifting through $\fbp[\fbp[E]]$ (respectively, $\fbl^{\uparrow p}[\fbl^{\uparrow p}[E]]$) (the proof covers $\fbl[E]$, but generalizes easily). This implies:

\begin{cor}\label{c:bap to pbap}
For a Banach space $E$, and $1 \leq p \leq \infty$, $\fbp[E]$ has the $\mu$-Approximation Property if and only if it has the $\mu$-Positive Approximation Property.
\end{cor}

In general, for every Banach lattice $X$, there exist a contraction $u : X \to \fbl[X]$, and a contractive lattice homomorphism $v : \fbl[X] \to X$, so that $I_X = vu$. We do not know whether we can always achieve a factorization with a positive $u$. By \cite{OiTTT1}, certain lattices $X$ do not lattice embed into any free lattice, so we cannot take $u$ to be a lattice homomorphism in the above factorization (no matter what $v$ is).

\section{Duals of free Banach lattices}\label{s:duality}

In this section we extract information about $\fbp[E]^*$ and $\fbl^{\uparrow p}[E]^*$  from the convex hull descriptions of $\ball(\fbp[E])$. More specifically, we estimate norms of linear combinations of atoms in these dual spaces. It is well known (see e.g.~\cite[p.~111]{AB}) that $z^* \in Z^*$ is an atom if and only if it implements a lattice homomorphism $Z \to \Real$. 

It is easy to see that, if $E$ is a Banach space, then an element of $\fbp[E]^*$ or $\fbl^{\uparrow p}[E]^*$ is an atom if and only if it is of the form $\widehat{x^*}$, for some $x^* \in E^*$. Indeed, any $x^* : E \to \Real$ has a unique lattice homomorphic extension, namely $\widehat{x^*}$. 
Conversely, suppose $\psi$ is an atom, that is, a lattice homomorphism from $\fbp[E]^*$ or $\fbl^{\uparrow p}[E]^*$ to $\Real$. Let $x^* = \phi|_E$, and note that both $\psi$ and $\widehat{x^*}$ are lattice homomorphisms extending $x^*$. However, $x^*$ has a unique extension to a free Banach lattice, hence $\psi = \widehat{x^*}$.

Thus, we shall estimate $\| \sum_{j=1}^N \widehat{x_j^*} \|$, with $x_1^*, \ldots, x_N^* \in E^*$.

% We know that the finite sums of atoms (that is, elements $\sum_j \widehat{x_j^*}$) are dense in $\fbp[E]^*$, and moreover, such elements of norm not exceeding one are dense in the unit ball of $\fbp[E]^*$. Indeed, for finite dimensional $E$ a simple discretization argument will do. When $E$ is infinite dimensional, we can reduce to finite dimensions, as shown \Cref{inf dim fbp}.

% However, for $\fbl^{\uparrow p}[E]^*$, we only have functional representation for finite dimensional $E$.

\subsection{Sums of atoms in $\fbl^{\uparrow p}[E]^*$}\label{ss:dual of fb p-conv}
We establish:

\begin{prop}\label{fbl p-conv dual}
 If $x_1^*, \ldots, x_N^* \in E^*$, then we have
 \begin{align*}
 \sup \Big( \sum_k \max_\pm \| \sum_{j \in S_k} \pm x_j^* \|^{p'} \Big)^{1/p'}
 &
 \leq \big\|\sum_j \widehat{x_j^*}\big\|_{\fbl^{\uparrow p}[E]^*}
 \\ &
 \leq K_p \sup \Big( \sum_k \max_\pm \| \sum_{j \in S_k} \pm x_j^* \|^{p'} \Big)^{1/p'} .
 \end{align*}
 Here the supremum runs over all representations of $\{1, \ldots, N\}$ as a union of disjoint sets $S_k$; $K_p$ is the constant from \Cref{hull upper est}.
\end{prop}

\begin{proof}
% Perturbing the elements $x_j^*$ slightly, we can and do assume they are not scalar multiples of each other.
%  By the results of \Cref{s:fbp hull}, we can compute the norm of $\psi = \sum_j \widehat{x_j^*}$ (up to a constant) by considering the action on elements $\varphi = \vee_k |\delta_{e_k}| \in \fbl^{\uparrow p}[E]$, with $\sum_k \|e_k\|^p \leq 1$. 
For the sake of brevity, let $\psi = \sum_j \widehat{x_j^*}$. By \Cref{hull upper est},
\begin{equation}
\sup_{\varphi \in S'} \langle \psi, \varphi \rangle \leq \|\psi\| \leq K_p \sup_{\varphi \in S'} \langle \psi, \varphi \rangle , \, \, {\textrm{where}} \, \,
S' = \big\{ \vee_k |\delta_{e_k}| : \sum_k \|e_k\|^p \leq 1 \big\} .
\label{norm up p}
\end{equation}

First we obtain an upper estimate on $\|\psi\|$. Fix $\varphi = \vee_k |\delta_{e_k}|  \in S'$.
 For each $j$ let $\sigma(j)$ be the smallest $k$ for which $\varphi(x_j^*) = |\langle x_j^* , e_k \rangle|$; let $S_k = \{j : \sigma(j) = k\}$. Then
 \begin{align*}
 \langle \psi, \varphi \rangle 
 &
 = \sum_k \sum_{j \in S_k} |\langle x_j^* , e_k \rangle| = \sum_k \max_\pm \sum_{j \in S_k} \pm \langle x_j^* , e_k \rangle \leq \sum_k \max_\pm \big\| \sum_{j \in S_k} \pm x_j^* \big\| \|e_k\| 
 \\
 &
  \leq \big(\sum_k \max_\pm \big\| \sum_{j \in S_k} \pm x_j^* \big\|^{p'} \Big)^{1/p'} \big( \sum_k \|e_k\|^p \big)^{1/p} \leq \big(\sum_k \max_\pm \big\| \sum_{j \in S_k} \pm x_j^* \big\|^{p'} \Big)^{1/p'} .
 \end{align*}
%  Taking the supremum o
By the right side of \eqref{norm up p},
 $$
 \|\psi\| \leq K_p \Big( \sum_k \| \max_\pm \sum_{j \in S_k} \pm x_j^* \|^{p'} \Big)^{1/p'} .
 $$
 
 The reverse inequality is established in a similar fashion. Suppose $(S_k)$ are disjoint subsets of $\{1, \ldots, N\}$, and
 $$
 \sum_k \big\| \sum_{j \in S_k} \pm x_j^* \big\|^{p'} = 1
 $$
 for some choice of $\pm$.  We shall show that, for $\psi = \sum_j \widehat{x_j^*}$, we have $\|\psi\| \geq 1$. To this end, fix $c \in (0,1)$; for each $k$, find a norm one $x_k \in E$ with $\sum_{j \in S_k} \pm \langle x_j^* , x_k \rangle > c \big\| \sum_{j \in S_k} \pm x_j^* \big\|$.
 Let $\alpha_k = \big\| \sum_{j \in S_k} \pm x_j^* \big\|^{p'/p}$, then $\sum_k \alpha_k^p = 1$, hence $\varphi = \vee_k \alpha_k |\delta_{x_k}|$ has norm not exceeding $1$. 
 Further, $ \langle \psi, \varphi \rangle = \sum_j \vee_i \alpha_i | \langle x_j^* , x_i \rangle|$. For each $j$, find the unique $k$ for which $j \in S_k$, then $\vee_i \alpha_i |\langle x_j^* , x_i \rangle| \geq \alpha_k |\langle x_j^* , x_k \rangle|$. 
 
 Recall that $p'/p + 1 = p'$, hence $\alpha_k \big\| \sum_{j \in S_k} \pm x_j^* \big\| = \big\| \sum_{j \in S_k} \pm x_j^* \big\|^{p'}$, and so,
 \begin{align*}
 \|\psi\| 
 &
 \geq \langle \psi, \varphi \rangle \geq \sum_k \alpha_k \sum_{j \in S_k} | \langle x_j^* , x_k \rangle | 
 \\
 &
 \geq c \sum_k \alpha_k \big\| \sum_{j \in S_k} \pm x_j^* \big\| = c \sum_k \big\| \sum_{j \in S_k} \pm x_j^* \big\|^{p'} = c ,  
 \end{align*}
 estimating $\|\psi\|$ from below ($c$ can be arbitrarily close to $1$).
\end{proof}

\subsection{Sums of atoms in $\fbp[E]^*$}\label{ss:dual of fbp}

% As before, fix $x_1^*, \ldots, x_N^* \in E^*$; our goal is to compute the norm of $\psi = \sum_j \widehat{x_j^*}$ in $\fbp[E]^*$. To this end some notation is needed. Denote by $(\sigma_j)_{j=1}^n$ the canonical basis of $\ell_p^n$; abusing the notation slightly, we shall use the same letter for different values of $n$ and $p$.

\begin{prop}\label{fbp dual}
For $1 \leq p \leq \infty$ and  $x_1^*, \ldots, x_N^* \in E^*$,
$$\big\|\sum_j \widehat{x_j^*}\big\|_{\fbp[E]^*} = \sup \Big\{ \|\sum_j u_j \otimes x_j^*\|_{\ell_{p'}^m(E^*)} : u_1, \ldots, u_N \in \ball(\ell_{p'}^m) \Big\} .$$
%  $$\|\psi\| = \sup \big\{ \|(A \otimes I_{E^*})(\sum_j \sigma_j \otimes x_j^*)\|_{\ell_{p'}^m(E^*)} : \| A : \ell_1^N \to \ell_{p'}^m\| \leq 1 \big\} . $$
\end{prop}

% Before proceeding to the proof, note that the operator $A$ has a matrix representation $(a_{ij})_{1 \leq i \leq m, 1 \leq j \leq N}$, via $A \sigma_j = \sum_i a_{ij} \sigma_i$. The condition $\| A : \ell_1^N \to \ell_{p'}^m\| \leq 1$ is equivalent to $\vee_j \sum_i |a_{ij}|^{p'} \leq 1$.

It may be instructive to consider \Cref{fbp dual} for ``extreme'' values of $p$.

$\bullet \,\, p = \infty$: $\big\|\sum_j \widehat{x_j^*}\big\| = \sum_j \|x_j^*\|$. Indeed, in the above notation, $\|\sum_j u_j \otimes x_j^*\| \leq \sum_j \|x_j^*\|$, and the equality is attained when $(u_j)$ form the canonical basis in $\ell_1^n$.

$\bullet \,\, p = 1$: % we are taking the supremum over all contractions $A : \ell_1^N \to \ell_\infty^m$. \Cref{fbp dual} tells us that
$\big\|\sum_j \widehat{x_j^*}\big\| = \sup_\pm \| \sum_j \pm x_j^*\| = \|(x_j^*)\|_{1,\textrm{weak}}$. % with the supremum running over all choices of sign. % sequences $(\alpha_j)$ with $\vee_j |\alpha_j| \leq 1$. % , or in other words, $\|\psi\| = 1$.

% we are taking the supremum over all contractions $A : \ell_1^N \to \ell_1^m$. As every Banach space is a subquotient of $\ell_1$, $\|(A \otimes I_{E^*})(\sum_j \sigma_j \otimes x_j^*)\| \leq \|\sum_j \sigma_j \otimes x_j^*\|$, hence, as $A$ can be taken to be the identity,
% $$ \|\psi\| = \|\sum_j \sigma_j \otimes x_j^*\|_{\ell_1^N(E^*)} = \sum_j \|x_j^*\| . $$
% This is not new: $\fbl^{(\infty)}[E]$ is an AM-space, and therefore, its dual has to be an AL-space.

\begin{proof}
% By approximation, we can and do assume that $x_j^*$'s are not scalar multiples of each other. 
For brevity let $\psi = \sum_j \widehat{x_j^*}$. By \Cref{hull}, it suffices to compute
 $$
 \sup \big\{ \big\langle \psi , (\sum_i |\delta_{e_i}|^p)^{1/p} \big\rangle : \sum_{i=1}^m \|e_i\|^p \leq 1 \big\} .
 $$
 However,
 \begin{align*}
 \big\langle \psi , (\sum_i |\delta_{e_i}|^p)^{1/p} \big\rangle 
 &
 = \sum_j \big( \sum_i |\langle x_j^* , e_i \rangle|^p \big)^{1/p} 
 \\
 &
 = \sup \big\{ \sum_j \sum_i a_{ij} \langle x_j^* , e_i \rangle : \vee_j \sum_i |a_{ij}|^{p'} \leq 1 \big\} ,
 \end{align*}
 and consequently,
 \begin{align*}
     \|\psi\| 
     & 
     = \sup \big\{ \sum_j \sum_i a_{ij} \langle x_j^* , e_i \rangle : \vee_j \sum_i |a_{ij}|^{p'} \leq 1, \sum_i \|e_i\|^p \leq 1 \big\}
     \\
     &
     = \sup \big\{ \sum_i \big\langle \sum_j a_{ij} x_j^* , e_i \big\rangle : \vee_j \sum_i |a_{ij}|^{p'} \leq 1, \sum_i \|e_i\|^p \leq 1 \big\}
     \\
     &
     = \sup \big\{ \big( \sum_i \big\|\sum_j a_{ij} x_j^*\big\|^{p'} \big)^{1/p'} : \vee_j \sum_i |a_{ij}|^{p'} \leq 1 \big\} ,
 \end{align*}
 and the right hand side equals
$$\sup \big\{ \|\sum_j u_j \otimes x_j^*\|_{\ell_{p'}^m(E^*)} : \max_j \|u_j\| \leq 1 \big\} $$
with $u_j = (a_{ij})_i$.
%  $$\sup \big\{ \|(A \otimes I_{E^*})(\sum_j \sigma_j \otimes x_j^*)\|_{\ell_{p'}^m(E^*)} : \| A : \ell_1^N \to \ell_{p'}^m\| \leq 1 \big\} . \qedhere $$
\end{proof}

The importance of finite linear combinations of atoms is made clear by the following result:

\begin{prop}\label{denseness}
 If $E$ is finite dimensional, then the elements of the form $\sum_{j=1}^N \widehat{x_j^*}$, of norm not exceeding $1$, are weak$^*$ dense in $\ball(\fbp[E]^*)_+$.
\end{prop}

\begin{proof}%[Proof of \Cref{denseness} for $\dim E < \infty$]
 Recall, from \Cref{bap}, that $\fbp[E]$ determines a convenient norm on $C(\sphere(E^*))$, with the controlling uniformly continuous family $\mathcal C$ (with modulus of continuity denoted by $\boldom$), consisting of functions $g = \big( \sum_i |\delta_{e_i}|^p \big)^{1/p}$, with $\sum_i \|e_i\|^p = 1$.
 Any element of $\ball(\fbp[E]^*)_+$ can be identified with a positive measure $\mu$, satisfying $\langle \mu, g \rangle \leq 1$, for any $g$ as above.
 
 Consider the family ${\mathcal A}$ of pointed partitions of $\sphere(E^*)$ into finitely many measurable subsets; that is, each $\alpha \in {\mathcal A}$ consists of a partition of $\sphere(E^*)$ into $S_1, \ldots, S_n$, and of a choice of $x_i^* \in S_i$. Let $\diam \alpha = \max_i \diam S_i$, and $\mu_\alpha = \sum_i \mu(S_i) \widehat{x_i^*}$. 
 Clearly $\mu_\alpha \to \mu$ in the weak$^*$ topology as $\diam \alpha \to 0$. It remains to show that 
 $$ \|\mu_\alpha\|_{\fbp[E]^*} \leq 1 + \dim E \cdot \boldom(\diam \alpha) . $$
 To this, note first that $\| \cdot \|_{\fbp[E]} \leq \dim E \cdot \| \cdot \|_\infty$, hence $\mu(\sphere(E^*)) \leq \dim E$. Further, for any $g \in {\mathcal{C}}$,
 \begin{align*}
 \big| \langle \mu, g \rangle - \langle \mu_\alpha, g \rangle \big| 
 &
 \leq \sum_{i=1}^n \int_{S_i} | g(e^*) - g(x_i^*) | \, d\mu 
 \\
 &
 \leq \sum_i \boldom(\diam S_i) \mu(S_i) \leq \dim E \cdot \boldom(\diam \alpha) ,
 \end{align*}
 which yields $\|\mu_\alpha\|_{\fbp[E]^*} \leq 1 + \dim E \cdot \boldom(\diam \alpha)$.
\end{proof}

% For general $E$, it suffices to show that the norm in $\fbp[E]$ can be computed ``locally.'' % We reproduce a reasoning of \cite{OiTT2}.

\section{Embeddings of free lattices}\label{s:embeddings}

In this section we show that the norm of any element of $\fbp[E]$ is determined in $\fbp[G]$, where $G$ is a (sufficiently large) finite dimensional subspace of $E$.

\begin{prop}\label{p:evaluate norm}
For any $f \in \fvl[E]$ and $\varepsilon > 0$, there exists a finite dimensional $G \subset E$ so that $\|f\|_{\fbp[E]} \leq (1+\varepsilon) \|f\|_{\fbp[G]}$.
\end{prop}

Begin the proof with:

\begin{lem}\label{l:extend operators}
 Suppose $1 \leq p \leq \infty$, and $F$ is a finite dimensional subspace of a Banach space $E$. Then for any $\varepsilon > 0$ there exist a finite dimensional $G$ $(F \subset G \subset E)$ with the following property: if $u \in B(F,L_p)$ is such that any extension $\widetilde{u} \in B(E,L_p)$ has norm at least $C$, then any extension $u' \in B(G,L_p)$ has norm at least $C/(1+\varepsilon)$.
\end{lem}

\begin{proof}
 Fix $\delta > 0$, and consider $u : F \to L_p$. By \cite[Lemma 2.2]{OiIP}, there exists an $M$ (depending on $\delta$ and $\dim F$ only) so that $L_p$ contains a sublattice $Z$, lattice isometric to $\ell_p^N$ (with $N \leq M$), so that there exists a contractive projection $Q$ from $L_p$ onto $Z$, and, for every $g \in u(F)$, $\dist(g, Z) \leq \delta \|g\|$. 
 Consequently, for any $f \in F$, $\dist(uf, Z) \leq \delta \|f\| \|u\|$. %By \cite[Lemma 2.2]{OiIP} again,
 As $\|I-Q\| \leq 2$, we conclude (again, as in \cite[Lemma 2.2]{OiIP}) that $\|(I-Q)u\| \leq 2 \delta \|u\|$. 
 As $F$ is $\sqrt{\dim F}$-complemented in any ambient space (Kadec-Snobar Theorem), $(I-Q)u$ extends to an operator from $E$ to $L_p$, of norm not exceeding $2 \sqrt{\dim F} \delta \|u\|$. Consequently, we have to show the following:  for any $\varepsilon > 0$ and $M \in \Nat$ there exist a finite dimensional $G$ with the following property: if $u \in B(F,\ell_p^M)$ is such that any extension $\widetilde{u} \in B(E,\ell_p^M)$ has norm at least $C$, then any extension $u' \in B(G,\ell_p^M)$ has norm at least $C/(1+\varepsilon)$.
 
 By a small perturbation argument, it suffices to show that $G$ as above exists for any finite collection of operators $u$. Since the span of finitely many finite dimensional spaces is again finite dimensional, it suffices to show the existence of a $G$ as above for a single operator $u$.
 
Suppose, for the sake of contradiction, that for any finite dimensional $G$ with $F \subset G \subset E$ there exists $u_G \in B(G,\ell_p^M)$ of norm less that $C'$. Taking an ultraproduct (as in \cite[Chapter 8]{DJT}), we achieve a contradiction by constructing $\widetilde{u} \in B(E,\ell_p^M)$, extending $u$, and of norm not exceeding $C'$.

The anonymous referee suggested an alternative conclusion to the proof, which avoids ultraproducts altogether, and instead resembles the Lindenstrauss compactness argument. 
Consider the set ${\mathcal{K}}$ of all functions $\ball(E) \to C' \ball(\ell_p^M)$, equipped with pointwise topology; by Tychonoff theorem, this set is compact.
Denote by ${\mathcal{G}}$ the set of finite dimensional spaces $G$ with $F \subset G \subset E$; denote by ${\mathcal{K}}_G$ the set of all elements of ${\mathcal{K}}$ which extend $u|_{\ball(E)}$, and whose restriction to $\ball(G)$ is linear. By the above reasoning, ${\mathcal{K}}_G \neq \emptyset$. 
By the finite intersection property, $\cap_{G \in \mathcal{G}} {\mathcal{K}}_G \neq \emptyset$; elements of this intersection correspond to the desired extensions $\widetilde{u}$.
\end{proof}

\begin{cor}\label{c:extend summing norms}
Suppose $F$ is a finite dimensional subspace of $E$, and $\varepsilon > 0$. Then there exists a finite dimensional $G$ with $F \subset G \subset E$ with the following property. If $N \in \Nat$, $(y_i^*) \in (F^*)^N$, and $z_1^*, \ldots, z_N^* \in G^*$ extend $y_1^*, \ldots, y_N^*$, then the same functionals also have an extensions $x_1^*, \ldots, x_N^* \in E^*$ with $\|(x_i^*)\|_{p,{\textrm{weak}}} \leq (1+\varepsilon) \|(z_i^*)\|_{p,{\textrm{weak}}}$.
\end{cor}

\begin{proof}[Sketch of a proof]
Identify $N$-tuples with operators into $\ell_{p'}$; for instance, $(y_i^*) \in F^N$ corresponds to the operator $T : F \to \ell_{p'} : y \mapsto ( \langle y_i^*, y \rangle )_i$, and $\|T\| = \|(y_i^*)\|_{p,{\textrm{weak}}}$.
\end{proof}

\begin{proof}[Proof of \Cref{p:evaluate norm}]
We can write $f$ as a lattice expression in $\delta_{x_1}, \ldots, \delta_{x_n}$, with $x_1, \ldots, x_n \in E$; let $F = \spn[x_1, \ldots, x_n]$. Note that, for each $x^* \in E^*$, $f(x^*)$ depends only on $x^*|_F$. Abusing the notation slightly, we can write $f(y^*)$ with $y^* \in F^*$.

Denote by $S$ the (convex) set of all $N$-tuples $y_1^*, \ldots, y_N^* \in F^*$ ($N \in \Nat$) for which there exist extensions $x_1^*, \ldots, x_N^* \in E^*$ so that $\|(x_j^*)\|_{p,{\mathrm{weak}}} < 1$. Then
$$
\|f\|_{\fbp[E]} = \sup \Big\{ \big( \sum_j \big| f(y_j^*) \big|^p \big)^{1/p} : (y_j^*) \in S \Big\} .
$$
By the above, there exists a finite dimensional $G$ with $F \subset G \subset E$, so that $S \subset S' \subset (1+\varepsilon) S$, where $S'$ is defined just like $S$, except for we require extensions from $F$ to $G$, not to $E$. Then
$$
\|f\|_{\fbp[G]} = \sup \Big\{ \big( \sum_j \big| f(y_j^*) \big|^p \big)^{1/p} : (y_j^*) \in S' \Big\} \leq (1+\varepsilon) \|f\|_{\fbp[E]} .
\qedhere
$$
\end{proof}

\begin{rem}
The preceding proof shows that the same $G$ works for all $f \in \fvl[E]$ % arising from
which are represented as lattice expressions in $\delta_{x_1}, \ldots, \delta_{x_N}$, with $x_1, \ldots, x_N \in F$. % in elements of a given finite dimensional $F$.
\end{rem}

\section{Extensions of linear operators to lattice homomorphisms}\label{s:maps}

If $E$ and $F$ are Banach spaces, then, as shown in \cite[Section 2]{JLTTT}, any $u \in B(E,F)$ extends, uniquely, to a lattice homomorphism $\overline{u} : \fvl[E] \to \fvl[F]$, sending $F(\delta_{x_1}, \ldots, \delta_{x_N})$ to $F(\delta_{ux_1}, \ldots, \delta_{ux_N})$ ($x_1, \ldots, x_N \in E$, and $F$ is a lattice expression).
This section aims to estimate the norm of $\overline{u}$ when $\fvl[E]$ ($\fvl[F]$) is viewed as a sublattice of $\fbl^{(q)}[E]$ or $\fbl^{\uparrow q}[E]$ (respectively,  of $\fbp[F]$ or $\fbl^{\uparrow p}[F]$). In other words, we want to determine when $u$ has a bounded extension $\overline{u} : \fbl^{(q)}[E] \to \fbp[F]$ (or similarly with free upper estimate lattices), and to estimate $\|\overline{u}\|$.

If $q \leq p$, then any $p$-convex Banach lattice is also $q$-convex, with the same constant (see e.g.~\cite[Proposition 1.d.5]{LT2}). By the defining property of $\fbl^{(q)}[E]$, then any $u : E \to F$ has a (necessarily unique) lattice homomorphic extension $\overline{u} : \fbl^{(q)}[E] \to \fbp[F]$, with $\|\overline{u}\| = \|u\|$.

Below, we focus on the interesting case of $p < q$. To proceed, recall that an operator $u \in B(X,Y)$ is called \emph{$(q,p)$-mixing} ($u \in M_{q,p}(X,Y)$) if there exists $C>0$ so that, for every regular Borel probability measure $\nu$ on $\ball(Y^*)$ (equipped with its weak$^*$ topology), there exists a regular Borel probability measure $\lambda$ on $\ball(X^*)$ (with weak$^*$-topology again), with the property that the inequality
$$
\Big( \int \big| \langle ux, y^* \rangle \big|^q \, d\nu(y^*) \Big)^{1/q} \leq C \Big( \int \big| \langle x, x^* \rangle \big|^p \, d\lambda(x^*) \Big)^{1/p} 
$$
holds for any $x \in X$. Define $\mu_{q,p}(u)$ to be the smallest $C$ as above. Importantly, $u \in M_{q,p}(X,Y)$ if and only if, for any $q$-summing $T : Y \to F$, $Tu$ is $p$-summing; in this case, $\mu_{q,p}(u) = \sup \pi_p(Tu)/\pi_q(T)$. For a detailed treatment of mixing operators, we refer the reader to \cite[Section 32]{DF}, \cite[Chapter 20]{Pietsch}, or to a neat summary in \cite{CD}.

The main result of this section is:

\begin{thm}\label{mixing}
 In the above notation, $u$ extends to a bounded lattice homomorphism $\overline{u} : \fbl^{(q)}[E] \to \fbp[F]$ if and only if $u^*$ is $(q,p)$-mixing. In this case, $\|\overline{u}\| = \mu_{q,p}(u^*)$.
\end{thm}

\begin{proof}
 (1) Suppose $\mu_{q,p}(u^*) < 1$, and show that, for any $\varphi \in \fvl[E]$ with $\|\varphi\|_q < 1$, we have $\|\overline{u} \varphi\|_p < 1$. This would show the boundedness of $\overline{u}$ on $\fvl[E]$, hence its extendability.
 
 By \Cref{p:fvl in fbp}, we can find $T : E^* \to Z$ (for some Banach space $Z$) so that $\pi_q(T) < 1$, and $|\varphi(e^*)| \leq \|T e^*\|$ for any $e^* \in E^*$. Then, for any $f^* \in F^*$,
 $$
 \big| [\overline{u} \varphi](f^*) \big| = \big| \varphi(u^* f^*) \big| \leq \|T u^* f^*\| .
 $$
 As $\pi_p(T u^*) < 1$, $\|\overline{u} \varphi\|_p < 1$, by \Cref{p:fvl in fbp}.
 
 (2) Now suppose $\|\overline{u}\| \leq 1$, and show that $\mu_{q,p}(u^*) \leq 1$. In other words, we have to show that, for any $q$-summing $T : E^* \to Z$, we have $\pi_p(T u^*) \leq \pi_q(T)$. Further, by Pietsch Factorization, it suffices to show that, whenever $\mu$ is a probability measure, and $S : E^* \to L_\infty(\Omega,\mu)$ is a contraction, then $\pi_p(JSu^*) \leq 1$ (here, $J : L_\infty(\mu) \to L_q(\mu)$ is the formal identity).
 
 Suppose $\FF$ is a partition of $\Omega$ into a disjoint union of sets $F_1, \ldots, F_N$; let $P_\FF$ be the corresponding conditional expectation. As $P_\FF$ is contractive, and commutes with $J$, we have $\pi_q(J P_\FF S) \leq 1$. On the other hand, we can order the set of partitions $\FF$ by inclusion; then $P_\FF \to id$ point-norm, hence $\lim_\FF J P_\FF S u^* = J S u^*$ point-norm, and therefore, $\pi_p(JSu^*) \leq \liminf_\FF \pi_p(J P_\FF S u^*)$. Thus, it suffices to show that $\pi_p(J P_\FF S u^*) \leq 1$ for any $\FF$.
 
 Henceforth, fix a partition $\FF$. Let $a_i = \mu(F_i)$. Denote by $W_\infty$ and $W_q$ the spans of $\chi_{F_i}$ ($1 \leq i \leq N$) in $L_\infty(\mu)$ and $L_q(\mu)$, respectively. We have isometries
 $$
 U_\infty : W_\infty \to \ell_\infty^N : \chi_{F_i} \mapsto \sigma_i \, {\textrm{  and  }} \,
 U_q : W_q \to \ell_q^N : \chi_{F_i} \mapsto a_i^{1/q} \sigma_i 
 $$
 (here $\sigma_1, \ldots, \sigma_N$ form the canonical bases in $\ell_\infty^N$ or $\ell_q^N$). Further, define $\Delta : \ell_\infty^N \to \ell_q^N : \sigma_i \mapsto a_i^{1/q} \sigma_i$ (recall that $\sum_i a_i = 1$). Then $J P_\FF S = U_q^{-1} \Delta U_\infty P_\FF S$. Therefore, it suffices to show that $\pi_p(T u^*) \leq 1$ for any $T = DR$, where $R : E^* \to \ell_\infty^N$ is a contraction, and $D : \ell_\infty^N \to \ell_q^N$ is a diagonal operator, with diagonal entries $(d_1, \ldots, d_N)$ satisfying $\sum_i d_i^q \leq 1$ (note that $\|D\| \leq 1$, hence also $\pi_q(T) \leq 1$).
 
 We can write $R e^* = ( \langle e_i^{**}, e^* \rangle )_i$, with $e_i^{**} \in \ball(E^{**})$.
 We can find a net $\mathfrak A$ so that, for each $i$ and $\alpha \in \mathfrak A$, we have $e_{i \alpha} \in \ball(E)$, and $e_{i \alpha} \to e_i^{**}$ in the weak$^*$ topology of $E^{**}$. The functionals $(e_{i\alpha})_{i=1}^N$ give rise to $R_\alpha : E^* \to \ell_\infty^N$, and $R_\alpha \to R$ point-norm. Consequently, $D R_\alpha u^* \to D R u^*$ point-norm, hence it suffices to show that $\pi_p(D R_\alpha u^*) \leq 1$ for every $\alpha$.
 
 Consider the function $\varphi_\alpha = \big( \sum_i d_i^q \big| \delta_{e_{i \alpha}} \big|^q \big)^{1/q} \in \fbp[E]$. Then $\|\varphi_\alpha\|_q \leq \big( \sum_i d_i^q \|e_{i\alpha}\|^q \big)^{1/q} \leq 1$. Thus, $\|\overline{u} \varphi_\alpha\|_p \leq 1$ -- that is, for any $f_1^*, \ldots, f_M^* \in F^*$ with $\|(f_j^*)\|_{p,{\textrm{weak}}} \leq 1$, we have $\sum_j \big|[\overline{u} \varphi_\alpha](f_j^*)\big|^p \leq 1$. However, for $f^* \in F^*$,
 $$
 [\overline{u} \psi_\alpha](f^*) = \varphi_\alpha (u^* f^*) = \big( \sum_i d_i^q | \langle u^* f^*, e_{i \alpha} \rangle |^q \big)^{1/q} = \|D R_\alpha u^* f^*\| ,
 $$
 hence for $(f_j^*)$ as above, we have $\sum_j \|D R_\alpha u^* f_j^*\|^p \leq 1$. As the latter inequality holds whenever $\|(f_j^*)\|_{p,{\textrm{weak}}} \leq 1$, we conclude that $\pi_p(D R_\alpha u^*) \leq 1$, as desired.
\end{proof}

As an application, we estimate $\|id : \fbl^{(q)}[E] \to \fbp[E]\|$, when $E$ is finite dimensional.

\begin{thm}\label{t:fin dim formal id}
 If $E$ is $n$-dimensional, then $\|id : \fbl^{(q)}[E] \to \fbp[E]\| \leq n^\alpha$, where $\alpha = \alpha(p,q)$ equals $0$ if $p \geq q$, $1/p - 1/q$ if $p \leq \min\{2,q\}$, and $p(1/p - 1/q)/2$ if $2 < p < q$.
\end{thm}

The above statement is obvious when $q \leq p$; for $q > p$, it is obtained by combining \Cref{mixing} with the following result, which may be of interest in its own right.

\begin{prop}\label{p:fin dim mixing}
 Suppose $X$ and $Y$ are Banach spaces, and $q > p$. Then, for any rank $n$ contraction $u : X \to Y$, we have $\mu_{q,p}(u) \leq n^\alpha$, with $\alpha$ as above.
\end{prop}

The proof of this proposition requires a change of density argument. Suppose $(\Omega,\mu)$ is a probability measure space, and $Z_q$ is a subset of $L_q(\Omega,\mu)$. Suppose $h$ is a non-negative function on $\Omega$, satisfying $\int h \, d\mu = 1$. Find a set $\Omega' \subset \Omega$ so that $h > 0$ $\mu$-a.e.~on $\Omega'$, and $h = 0$ $\mu$-a.e.~on $\Omega \backslash \Omega'$. Then $\mu' = h \mu$ is a probability measure on $\Omega'$.
Assume that all functions from $Z_q$ vanish outside of $\Omega'$. Then the map $\Phi_{qh} : Z_q \to L_q(\mu') : f \mapsto h^{-1/q} f$ is an isometry from $Z_q \subset L_q(\mu)$ onto $Z_q' \subset L_q(\mu')$.

For a given $h$, let $\CC(h,Z_q,p,q) = \sup_{f \in Z_q}\|\Phi_{qh} f\|_q/\|\Phi_{qh} f\|_p$, $\CC(Z_q, q,p) = \inf_h \CC(h,Z_q,p,q)$. For $n \in \Nat$, let $\CO(n, q,p) = \sup_{\dim Z_q = n} \CC(Z_q,p,q)$ (note that, when $p$ and $q$ are fixed, then $\CO(n, q,p)$ is a non-decreasing function of $n$).

\Cref{p:fin dim mixing} would then follow by combining two lemmas below.

\begin{lem}\label{l:mixing vs density}
 For $u$, $p$, $q$ as in \Cref{p:fin dim mixing}, $\mu_{q,p}(u) \leq \CO(n, q,p)$.
\end{lem}

\begin{lem}\label{l:density est}
 For $p$, $q$, $\alpha$ as in \Cref{p:fin dim mixing}, $\CO(n, q,p) \leq n^\alpha$.
\end{lem}

\begin{proof}[Proof of \Cref{l:mixing vs density}]
 By Pietsch Factorization, it suffices to show that $\pi_p(Tu) \leq \CO(n, q,p) \pi_q(T)$ for any operator $T = JR$, where $R : Y \to L_\infty(\mu)$ is a contraction, $\mu$ is a probability measure, and $J : L_\infty(\mu) \to L_q(\mu)$ is the formal identity.
 
 Let $Z_q = JRu(X) \subset L_q(\mu)$, and note that $m = \dim Z_q \leq n$. Fix $C' >  \CO(n, q,p)$, and find $h$ so that $\CC(h,Z_q,q,p) < C'$. As before, let $Z_q' = \Phi_{qh}(Z_q) \subset L_q(\mu')$; denote by $Z_p'$ the same subspace, equipped with the norm inherited from $L_p(\mu')$. 
 Denote the formal identity $Z_q' \to Z_p'$ by $j$. Then $j$ is contractive (due to $\mu'$ being a probability measure), and $\|j^{-1}\| < C'$. Further, denote by $W$ the restriction of $\Phi_{qh}$ to an isometry $Z_q \to Z_q'$.
 Also, denote by $J'$ the formal identity from $L_\infty(\mu')$ to $L_q(\mu')$, and by $V$ the natural embedding of $L_\infty(\mu)$ to $L_\infty(\mu')$. Then, for any essentially bounded function $g$ on $\Omega$ which vanishes outside of the support of $h$, we have $\Phi_{qh} J g = J' V g$. Let $S = V R u$ (this is well defined, since $h$ does not vanish on the image of $Ru$).
 
 Use the notation of the preceding paragraph to factor $Tu$:
 $$
 Tu = W^{-1} j^{-1} \circ j \Phi_{qh} J \circ Ru = W^{-1} j^{-1} \circ j J' \big|_{\ran S} \circ S
 $$
 Note that $j J' \big|_{\ran S}$ is a restriction of the formal identity $L_\infty(\mu') \to L_p(\mu')$, hence the right side exhibits a $p$-summing Pietsch decomposition of $Tu$. Therefore, $\pi_p(Tu) \leq \|W^{-1}\| \|j^{-1}\| \|S\| < C'$. To finish the proof, recall that $C'$ can be arbitrarily close to $\CO(n, q,p)$.
\end{proof}

\begin{proof}[Proof of \Cref{l:density est}]
(1) $q > 2 \geq p$: we shall obtain a desirable change of density by applying Lewis's Lemma (see \cite{Lew}, or a survey \cite{JL}), in the same manner as it is done in e.g.~\cite[Lemma 2]{Schechtman87}.
Suppose $Z_q$ is an $n$-dimensional subspace of $L_q(\mu)$. Use \cite{Lew} to find a basis $(x_i)_{i=1}^n \subset Z_q$ so that, for any $(a_i)$, we have
$$
\frac1n \sum_i a_i^2 = \int \Big| \sum_i a_i x_i \Big|^2 h^{q-2} \, d\mu , \, \, {\textrm{where}} \, \, h = \big( \sum_i |x_i|^2 \big)^{1/2} .
$$
Note that, for any sequence $a_i = \pm 1$, we have
$$
1 = \int \big( \sum_{i,j} a_i a_j x_i x_j \big) h^{q-2} \, d\mu ;
$$
averaging over all such $\pm 1$-strings, we conclude that $\int h^q \, d\mu = 1$.

Let $d \mu' = h^q \, d\mu$ (by the above, this is a probability measure), $x_i' = x_i/h \in L_q(\mu')$, and consider a change of density $f \mapsto f/h$. Then $\sum_i |x_i'|^2 = 1$, and, for any $(a_i)$,
$$
\Big\| \sum_i a_i x_i' \Big\|_{L_2(\mu')}^2 = \int \Big| \sum_i a_i x_i' \Big|^2 h^q \, d\mu =
\int \Big| \sum_i a_i x_i \Big|^2 h^{q-2} \, d\mu = \frac1n \sum_i a_i^2 .
$$
Therefore,
\begin{align*}
 \Big\| \sum_i a_i x_i' \Big\|_{L_q(\mu')}^q
 & = \int \Big( \sum_i a_i x_i' \Big)^q \, d \mu' = \int \Big( \sum_i a_i x_i' \Big)^2 \Big( \sum_i a_i x_i' \Big)^{q-2} \, d \mu'
\\ &
 \leq \int \Big( \sum_i a_i x_i' \Big)^2 \Big( \sum_i a_i^2 \Big)^{(q-2)/2} \Big( \sum_i x_i'^2 \Big)^{(q-2)/2} \, d \mu '
\\ &
 = \Big( n \Big\| \sum_i a_i x_i' \Big\|_{L_2(\mu')}^2 \Big)^{(q-2)/2} \int \Big( \sum_i a_i x_i' \Big)^2 \, d \mu' 
 \\ &
 = n^{(q-2)/2} \Big\| \sum_i a_i x_i' \Big\|_{L_2(\mu')}^q .
\end{align*}
Consequently, $\|f\|_q \leq n^\beta \|f\|_2$, for any $f \in Z_q'$ ($\beta = 1/2 - 1/q$). By H\"older Inequality, 
$$
\|f\|_2 \leq \|f\|_q^{1-\theta} \|f\|_p^\theta \, {\textrm{  with  }} \, \frac12 = \frac{1-\theta}q + \frac\theta{p} .
$$
Thus, $n^{-\beta} \|f\|_q \leq \|f\|_q^{1-\theta} \|f\|_p^\theta$, which leads to $\|f\|_q \leq n^{\beta/\theta} \|f\|_p = n^\alpha \|f\|_p$.

(2) $2 \geq q > p$: Again we use Lewis's Lemma, changing the density in such a way that $\|f\|_2 \leq n^\beta \|f\|_p$ (with $\beta = 1/p - 1/2$) holds for $f \in Z_q'$. We then deploy a similar interpolation argument:
$$
\|f\|_q \leq \|f\|_p^{1-\theta} \|f\|_2^\theta \, {\textrm{  with  }} \, \frac1q = \frac{1-\theta}p + \frac\theta{2} ,
$$
hence $\|f\|_q \leq n^{\beta\theta} \|f\|_p$, and $\beta\theta = \alpha$.

(3) $q > p \geq 2$: instead of Lewis's Lemma, we use \cite[Lemma 2]{Schechtman87}. %, and its proof.
A change of density guarantees that $\|f\|_\infty \leq n^{1/2} \|f\|_2$, for any $f \in Z_2'$. By H\"older Inequality, 
$$
n^{(1-p/q)/2} \|f\|_p \geq \|f\|_p^{p/q} \|f\|_\infty^{1-p/q} \geq \|f\|_q . \qedhere
% \|f\|_q \leq \|f\|_p^{1-\theta} \|f\|_\infty^\theta \, {\textrm{  with  }} \, \frac1q = \frac{1-\theta}p + \frac\theta\infty .
$$
% So $\|f\|_q \leq n^{\theta/p} \|f\|_p$, and $\theta/p = \alpha$.
\end{proof}

\begin{rem}\label{r:optimality}
 \Cref{r:Banach Mazur} below shows that, in general, the estimate of \Cref{l:density est} is optimal.
\end{rem}

We can conclude that, for a finite dimensional $E$, $\fbp[E]$ and $\fbl^{\uparrow p}[E]$ are ``not too different.'' Note first that $id : \fbl^{\uparrow p}[E] \to \fbp[E]$ is contractive. For converse, we have:

\begin{cor}\label{fbp to upper p}
There exists a universal constant $C$ so that $\|id : \fbp[E] \to \fbl^{\uparrow p}[E]\| \leq C (\ln n)^{1/p}$, whenever $1 < p < \infty$, and $E$ is $n$-dimensional, with $n$ large enough.
\end{cor}

\begin{proof}
%  Fix $r \in [1,p)$ (to be specified later). 
Let $r = p(1- 1/\ln n)$, and factor $id : \fbp[E] \to \fbl^{\uparrow p}[E]$ as a composition $vu$, where $u = id : \fbp[E] \to \fbl^{(r)}[E]$ and $v = id : \fbl^{(r)}[E] \to \fbl^{\uparrow p}[E]$. 

First consider $p \in (1,2]$. By \Cref{t:fin dim formal id}, $\|u\| \leq n^{1/r-1/p}$.
 By definition, $\fbl^{\uparrow p}[E]$ has upper $p$-estimate with constant $1$, hence, by \cite[Theorem 1.f.7]{LT2}, it is $r$-convex with constant $K \leq K_1 (p/(p-r))^{1/r}$, where $K_1$ is a certain universal constant. Thus (see e.g.~\cite[Theorem 2.1]{OiTTT1}), $\|v\| \leq K$.
 Consequently,
 $$
 \|id : \fbp[E] \to \fbl^{\uparrow p}[E]\| \leq \|u\| \|v\| \leq K_1 \Big( \frac{p}{p-r} \Big)^{1/r} n^{1/r - 1/p} ,
 $$
%  Now, selecting $r = p(1- 1/\ln n)$ yields the desired result.
which yields the desired result.

For $p \in (2,\infty)$, \Cref{t:fin dim formal id} gives $\|u\| \leq n^{(1-r/p)/2}$, hence
$$
 \|id : \fbp[E] \to \fbl^{\uparrow p}[E]\| \leq \|u\| \|v\| \leq K_1 \Big( \frac{p}{p-r} \Big)^{1/r} n^{(1-r/p)/2} ;
 $$
 again, the conclusion of the corollary holds.
\end{proof}

\begin{rem}\label{r:Banach Mazur}
 Similarly to the well-known Banach-Mazur distance between Banach spaces (see e.g.~\cite{TJ}), one can define the \emph{lattice Banach-Mazur distance}, denoted by $d_\ell$. For lattices $X$ and $Y$, we set $d_\ell(X,Y) = \inf \|T\| \|T^{-1}\|$, where the infimum runs over all lattice isomorphisms $T : X \to Y$ (if $X$ is not lattice isomorphic to $Y$, we set $d_\ell(X,Y) = \infty$).
 \Cref{t:fin dim formal id} shows that, if $\dim E = n$, and $1 \leq \min\{p,q\} \leq 2$, then $d_\ell(\fbp[E], \fbl^{(q)}[E]) \leq n^{|1/p-1/q|}$.
 
 In general, this estimate is optimal. For instance, we claim that, for $E = \ell_1^n$, we have $d_\ell(\fbl^{(\infty)}[E], \fbl[E]) \geq n$.
 To establish this inequality, consider a contractive lattice homomorphism $T : \fbl[E] \to \fbl^{(\infty)}[E]$. Then $T^* : \fbl[E]^* \to \fbl^{(\infty)}[E]^*$ is interval-preserving, hence maps atoms to atoms. As in \cite{OiTTT1}, we observe that there exists a continuous positively homogeneous function $\Psi : E^* \to E^*$ so that $T^* \widehat{x^*} = \widehat{\Psi x^*}$, for any $x^* \in E^*$.
 
 Let $x_1^*, \ldots, x_n^*$ be the canonical basis of $E^* = \ell_\infty^n$; for $1 \leq i \leq n$ let $y_i^* = \Psi(x_i^*)$. Further, let $\psi = \sum_{i=1}^n \widehat{x_i^*}$.
 By the discussion following \Cref{fbp dual}, $\|\psi\|_{\fbl[E]^*} = \max_\pm \|\sum_i \pm x_i^*\| = 1$. On the other hand,
 $\|T^* \psi\|_{\fbl^{(\infty)}[E]^*} = \sum_i \|y_i^*\| \leq 1$.
 We conclude that, for some $i$, $\|y_i^*\| \leq 1/n$, hence
 $$
 \|T^{-1}\| = \|T^{*-1}\| \geq \frac{\|\widehat{x_i^*}\|}{\|T^*\widehat{x_i^*}\|} = \frac{\|x_i^*\|}{\|y_i^*\|} \geq n .
 $$
\end{rem}

\begin{rem}\label{r:optimality of log}
\Cref{fbp to upper p} shows that $d_\ell(\fbp[E], \fbl^{\uparrow p}[E]) \prec (\ln n)^{1/p}$. This dependence on $n$ cannot, in general, be improved. Indeed, fix $p \in (1,\infty)$. Let $E = \ell_{p\infty}^n$, with its Banach space norm:
 $$
 \| (\alpha_i)_{i=1}^n \| = \| (\alpha_i)_{i=1}^n \|_{p,\infty} = \sup_{S \subset \{1, \ldots, n\}} \frac1{|S|^{1/p'}} \sum_{i \in S} |\alpha_i| 
 $$
 (this is \eqref{eq:p infty new norm}, applied to the counting measure on $\{1, \ldots, n\}$).
 We shall show that the $p$-convexity constant of $\fbl^{\uparrow p}[E]$ is at least $c  (\ln n)^{1/p}$, for some constant (which may depend on $p$, but not on $n$). We do so by exhibiting $x_1, \ldots x_n \in E$ so that $\|x_k\| \prec 1$ for any $k \in \{1, \ldots, n\}$ (and consequently, $\big( \sum_k \|x_k\|^p \big)^{1/p} \prec n^{1/p}$), yet
 \begin{equation}
  \label{eq:norm est} 
  \Big\| \Big( \sum_k \big| \delta_{x_k} \big|^p \Big)^{1/p} \Big\| \succ (\ln n)^{1/p} n^{1/p}
 \end{equation}
% To this end, define, for $1 \leq j \leq n$, $\beta_j = j^{-1/p}$. 
For $1 \leq k \leq n$, we let $x_k = (\beta_{ik})_{i=1}^n$, where $\beta_{ik} = (i+k-1 \, \mod n)^{-1/p}$ (with the convention that $m \mod n \in \{1, \ldots, n\}$).
% $$ x_k = ( \beta_k, \beta_{k+1}, \ldots, \beta_n, \beta_1, \ldots, \beta_{k-1} ) . $$
% Clearly $\|x_k\| \sim 1$, for any $k$.
Note that, for $1 \leq j \leq n$, exactly $j$ coordinates of $x_k$ are greater or equal than $j^{-1/p}$, hence, by \eqref{eq:p infty norm},
$$
\triple{x_k}_{p,\infty} = \max_{1 \leq j \leq n} j^{-1/p} \cdot j^{1/p} = 1 .
$$
By the equivalence of the norms $\| \cdot \|_{p,\infty}$ and $\triple{ \cdot }_{p,\infty}$ (discussed in the paragraph preceding the proof of \Cref{hull upper est}), $\|x_k\| \sim 1$.

The standard theory of Lorentz spaces (see e.g. \cite{Dilworth}) tells us that $E^* = \ell_{p'1}^n$, with the norm
 $$
 \| (\alpha_i)_{i=1}^n \| = \sum_{i=1}^n \big(i^{1/p'} - (i-1)^{1/p'}\big) \alpha_i^* ,
 $$
 where $(\alpha_i^*)$ is the decreasing rearrangement of $(\alpha_i)$. 
Let $(e_i^*)_{i=1}^n$ be the canonical basis of $E^*$. For any $S \subset \{1, \ldots, n\}$, and any choice of $\pm$, we have $\|\sum_{i \in S} \pm e_i^*\| = |S|^{1/p'}$, hence, by \Cref{fbl p-conv dual}, $\|\sum_{i=1}^n \widehat{e_i^*}\|_{\fbl^{\uparrow p}[E]^*} \sim n^{1/p'}$. Consequently,
\begin{align*}
\Big\| \Big( \sum_k \big| \delta_{x_k} \big|^p \Big)^{1/p} \Big\| 
&
\geq \big\|\sum_{i=1}^n \widehat{e_i^*}\big\|^{-1} \Big\langle \Big( \sum_k \big| \delta_{x_k} \big|^p \Big)^{1/p} , \sum_i \widehat{e_i^*} \Big\rangle 
\\
&
\succ n^{-1/p'} \sum_i \big( \sum_k | \langle e_i^*, x_k \rangle |^p \big)^{1/p} = n^{-1/p'} \sum_i \big( \sum_k \beta_{ik}^p \big)^{1/p} 
\\
&
= n^{-1/p'} \sum_i \big( \sum_j j^{-1} \big)^{1/p} \sim n^{-1/p'} \cdot n (\ln n)^{1/p} ,
\end{align*}
thus establishing \eqref{eq:norm est}.
\end{rem}

\end{document}